\documentclass[12pt]{amsart}
\usepackage{amsmath}
\usepackage{amssymb}
\usepackage{amsfonts}
\usepackage{amsthm}
\usepackage{verbatim}
\usepackage{amscd}
\usepackage{cite}
\usepackage{leftidx}
\usepackage{enumerate}
\usepackage{txfonts}
\usepackage{manfnt}
\usepackage{amscd}
\usepackage[mathscr]{eucal}
\usepackage{hyperref}
\usepackage{datetime2}
\usepackage{mathpazo}
\def\cP{\mathcal P}
\def\cX{\mathcal X}
\def\cZ{\mathcal Z}
\def\cC{\mathcal C}

\newtheorem{thm}{Theorem} % [section]

\newtheorem*{thm*}{Theorem}
\newtheorem*{prop*}{Proposition}
\newtheorem{cor}[thm]{Corollary}
\newtheorem*{cor*}{Corollary}

\newtheorem{lem}[thm]{Lemma}
\newtheorem*{lem*}{Lemma}

\newtheorem*{claim*}{Claim}
\newtheorem{prop}[thm]{Proposition}

\theoremstyle{remark}

\newtheorem{rem}[thm]{Remark}
\newtheorem*{rem*}{Remark}
\newtheorem{crit-rem}[thm]{Critical remark}

\newtheorem{example}[thm]{Example}
\newtheorem*{example*}{Example}

\newtheorem*{defn*}{Definition}

\def\inv{^{-1}}

\def\refp #1.{(\ref{#1})}
\newcommand\carets [1]{\langle #1 \rangle}

\newcommand{\A}{\mathcal{A}}

\newcommand{\kk}{\mathbf{k}}

\newcommand{\ul}[1]{\underline {#1}}

\def\sbr #1.{^{[#1]}}
\def\sfl #1.{^{\lfloor #1\rfloor}}

\def\inv{^{-1}}
\def\?{{\bf{??}}}

\def\A{\Bbb A}

\def\C{\mathbb C}
\def\P{\mathbb P}
\def\N{\mathbb N}
\def\R{\mathbb R}
\def\Z{\mathbb Z}

\def\O{\mathcal O}

\def\rk{\text{rk}}

\def\g{\mathfrak g}

\def\1/2{\frac{1}{2}}

\def\I{\mathcal{ I}}

\def\2{{[2]}}
\def\l{\ell}
\def\nl{\newline}

\def\<{\langle}
\def\>{\rangle}

\def\2{{[2]}}
\def\l{\ell}

\def\scl #1.{^{\lceil#1\rceil}}
\def\spr #1.{^{(#1)}}
\def\sbc #1.{^{\{#1\}}}

\def\subpr#1.{_{(#1)}}

\def\beq{\begin{equation*}}
\def\eeq{\end{equation*}}

\def\g3{{\Gamma\spr 3.}}
\def\gg{{\Gamma\spr 2.}}

\def\gg{\mathbb G}

\newcommand{\eqspl}[2]{
\begin{equation}\label{#1}
\begin{split}
#2\end{split}\end{equation}}

\newcommand{\exseq}[3]{
0\to #1\to #2\to #3\to 0
}

\newcommand{\beginalphaenum}{
\begin{enumerate}\renewcommand{\labelenumi}{ }
\item \begin{enumerate}
}

\def\eex{\end{rm}\end{example}}
\begin{document} 

%\end{document}
%\title{Differential complexes, degeneracy\\ and Hodge theory on Poisson manifolds}
\title{Balanced curves and minimal rational connectedness\\  
	 on Fano hypersurfaces}
\author %{author}
{Ziv Ran}
%Partially supported by NSA Grant MDA904-02-1-0094} }
\thanks{arxiv.org/math.ag/2008.01235 }
%\date {\today}% \enddate
\date{\DTMnow}

%\affil University of California, Riverside\endaffil

\address {\nl UC Math Dept. \nl
	Skye Surge Facility, Aberdeen-Inverness Road
	\nl
	Riverside CA 92521 US\nl 
	ziv.ran @  ucr.edu\nl
	\url{https://profiles.ucr.edu/app/home/profile/zivran}
}

%\email {ziv.ran @ucr.edu}
\subjclass[2010]{14n25, 14j45, 14m22}

\begin{abstract}
	
	On a general  hypersurface $X$ of degree $n$ (resp. $<n$) in $\P^n$
	and any $e\geq n-1$ (resp. arbitrarily large $e$),
	we construct families of rational curves of degree $e$
	going through the maximal number of general points. This
	solves in many cases the problem of rational curve interpolation on $X$.
%we 	determine for infinitely many $k$ the minimal degree $e$ of a 
%	rational curve through a 
%	general  collection of $k$ points.  
%	In the case of a hypersurface of  index 1, our results hold
%	for all $k\geq 1$. In an appendix, M.C. Chang shows that
%	in the case of index $>1$, 
%	the density of the set of curve degrees $e$ covered by our method
%	is approximately $\frac{(n-d)(d-\frac{5}{2})}{(n-2)d}$.
	
%	We study vector bundles on certain rational trees and their smoothings.
%	We use this to construct rational curves with balanced normal bundle 
%	for infinitely many degrees 
%	$e\geq n-1$ on a general hypersurface of degree $d\leq n$ in $\P^n, n\geq 4$.
%	In case $d=n$, all degrees $e$ are covered. This has consequences as
%	to the lowest degree of a rational curve through general points and
%	to the separable higher rational connectedness of these hypersurfaces.
%	The proof is based on fan degeneration and its generalization.
\end{abstract}
\maketitle
As is well known, $q$ general points in the plane lie on a rational
curve of degree $\lceil\frac{q+1}{3}\rceil$ and none lower.
More generally, 
given a (polarized, rationally connected) variety $X$,
one is interested in the question of \emph{rational curve interpolation}
(aka minimal rational connectivity) on $X$: 
what is the minimal degree
of a rational curve on $X$
through a general collection of $q$ points? 
At least when $X$ is anticanonically polarized, there
is an obvious expected answer: namely, 
the minimal anticanonical degree 
of such a curve $C$ is the smallest such that
\[\chi(N_{C/X})=C.(-K_X)+\dim(X)-3\geq q(\dim(X)-1),\] i.e. it equals
\[\min\{C.(-K_X):C.(-K_X)\geq (q-1)(\dim(X)-1)+2\}.\] 
When the expected degree equals the actual degree we say that
$X$ is minimally rationally $(q-1)$-connected.
%In other words,  the expected
%minimal
%anticanonical degree is
%\[C.(-K_X)=(k-1)(\dim(X)-1)+2.\]
The results of this paper show (in any characteristic) that
 a general Fano hypersurface $X$ in 
projective space is (separably) minimally $(q-1)$- rationally connected
for infinitely many values of
 $q$ and even \emph{any} $q\geq 1$ if $X$
has index 1; moreover, in these cases the locus of rational curves through $q$ 
general points  is reduced of the
expected dimension. We proceed to describe the results precisely.

%Rational curves on (Fano) hypersurfaces $X$ of degree $d\leq n$ in $\P^n$
%have been much studied in recent years, especially by Joe Harris and his %students,
%see e.g. \cite{harris-roth-starr}, %\cite{coskun-riedl}\cite{kollar-rat-curves},
%\cite{riedl-yang}. 

Here $\O_C(k)$ denotes the line bundle of degree $k$ and $a\O_C(k)$ denotes
$\bigoplus\limits_1^a\O_C(k)$.
A rational curve $C\to X$ is said to be \emph{balanced}
if its normal bundle $N_C$ is a balanced bundle on $C\simeq\P^1$, which means,
setting $n=\dim(X)+1=\rk(N_C)+2$, that
 \[N_C\simeq (n-2-r^-)\O_C(a^-+1)\oplus r^-\O_C(a^-)\ \ \mathrm{for\  some}\ \  r^-> 0, a^-\in\Z.\]
 In that case
  \[a^-=\lfloor\deg(N_C)/(n-2)\rfloor=\lfloor (C.(-K_X)-2)/(n-2)\rfloor.\] 
 Geometrically, 
balancedness implies
that $C$ is movable to go through the expected- i.e.
maximal- number, viz. 
$(a^-+1)$, of general
points on $X$, hence 
balancedness is closely related to rational connectedness and its
generalizations. A (polarized) variety $X$ is 
{rationally $(q-1)$-connected}, $q\geq 1$,
if  there is a family of rational curves $\cC/B\to X$
such that the induced map $\cC^{q}/B\to X^{q}$
is dominant. $X$ is 
\emph{rationally
$(q-1,e)$-connected} if the curves can be taken to have polarized degree $e$.
The adjective 'separable' may be
 added to these properties if the induced map $\cC^{q}/B\to X^{q}$ 
 is separable (as well as dominant). This makes most sense if
 the polarization is $H=(-1/k_X)K_X$ for some \emph{index} $k_X\in\N$, where
 $k_X=n+1-d$ for a hypersurface $X$ of degree $d\leq n$ in $\P^n$.
Thus, 
the existence of a balanced rational curve
of degree $e$ on $X$ 
 is equivalent to separable
 $(a^-, e)$  -
rational
connectedness, where $a^-=\lfloor\frac{k_Xe-2}{n-2}\rfloor$.
Fixing $n, d$, we say that $e$ is \emph{point-minimal} when
\[  \frac{(n+1-d)(e-1)-2}{n-2}<\lfloor\frac{(n+1-d)e-2}{n-2}\rfloor,\eqno{(*)}
\]%\lfloor\frac{k_Xe-2}{n-2}\rfloor$, 
Whenever $e$ is point-minimal, $(a^-,e)$ rational connectedness 
  implies that $e$ is the minimal, as well as the expected,
degree of a rational curve through $a^-+1$ general points, 
in which case we say $X$ is \emph{minimally rationally $a^-$-connected}
and that the point-degree $a^-+1$ is \emph{interpolating} (for $X$). 
\par
Rational connectedness and $q$-connectedness, not necessarily minimal,
 of all Fano manifolds has been known since the 90s
(see Koll\'ar's book \cite{kollar-rat-curves}).
For general hypersurfaces $X$ of degree $d\leq n$ in $\P^n$,
 Chen and Zhu \cite{chen-zhu} and Tian \cite{tianz}
have proven that $X$ is separably rationally connected.
Some more precise results on existence of low-degree balanced rational curves 
(and consequently, minimal rational connectedness in low degrees)
for such hypersurfaces are given in  
\cite{coskun-riedl} and\cite{hypersurf}
(see also \cite{alzati-re-normal},
  \cite{shen-normal}). 
 
  In this paper 
we extend these results to cases of  curves of high degree $e$.
%It is convenient to set, for fixed $d\leq n$, 
%\[q_{\max}(e)=\lfloor\frac{e(n+1-d)-2}{n-2}\rfloor+1\] 
%\[e_{\min}(q)=\lceil\frac{(q-1)(n-2)+2}{n+1-d}\rceil.\]
%Some trivial remarks:
%\begin{itemize}\item
%$q_{\max}(e)$ is the maximum  (expected) number of general points
%that may lie on a rational curve $C$ of degree $e$ on a hypersurface $X$ of degree $d$ in $\P^n$,
%as well as the actual number of general points on a general $C$ if $C$ is balanced.
%\item Similarly $e_{\min}(q)$ is the minimum (expected) degree
%of a rational curve that may go through $q$ general points of $X$
%and the expected degree equals the actual degree iff $X$ is minimally $q$-connected.
%\item The two functions above are 'almost' mutual inverses:
% If $d\geq 3$ then $\frac{n+1-d}{n-2}<1$ so by elementary arithmetic
% $q_{\max}(e_{\min}(q))=q$
%while $e_{\min}(q_{\max}(e))\leq e$, with equality iff $q_{\max}(e-1)<q_{\max}(e)$.
%\item Thus, if $X$ contains a balanced rational curve of degree $e=e_{\min}(q)$
%then $X$ is minimally rationally $q$-connected (for general $e$ there may exist rational curves of degree $<e$
%through $q_{\max}(e)$ general points).\end{itemize}
%and to note that 
%if $d\geq 3$ then $q_{\max}(e)-1\leq q_{\max}(e-1)\leq q_{\max}(e)$.
Our main results are as follows (here 'general' refers to an open set,
depending on the values of the numerical parameters
involved, in the space of hypersurfaces,
 and having the results valid  for all values simultaneously requires 'very general'):
\begin{itemize}
	\item
For $d=n\geq 4$, $e\geq n-1$, a general $X$ contains balanced 
rational curves  of degree $e$  
 (Theorem \ref{d=n-thm}), hence a very general $X$ is separably minimally rationally
$q$-connected for all $q$  (Corollary \ref{ratcon-d=n-cor}). 
\item 
For each $3\leq d<n$ and $e, q$ satisfying certain arithmetical conditions,
%arbitrarily large
%$e$ depending on $d,n$, 
a general $X$
contains  balanced rational curves of degree $e$ and is
separably minimally rationally $q$-connected   (Theorem \ref{d<n-thm}). 
%consequently, for arbitrarily large $q$
%a general $X$ is separably rationally $q$-connected; 
%these results come with a lower estimate on the density of  
%set of $e$ and $q$ involved.
%hence if $e$ is point-minimal,
%$X$ is minimally separably $a^-$-rationally connected.
\par For example (cf. Example \ref{d=n-1-example}), 
if $d=n-1$ then a general $X$ is minimally $((n-1)k+1)$- (resp. $(n-1)k+2$)- rationally connected
for all $k\geq 1$ if $n$ is even (resp. odd) and $\geq 6$.
%*******$q:=q_{\max}(e)>q_{\max}(e-1)$,
% $X$ is separably minimally $q$- 
%rationally connected,  (???with
 %$q\nearrow\infty$ (Corollary \ref{d<n-cor})???).*****
More generally,   a result of M. C. Chang
 (see the Appendix) shows, for all $2<d<n-1$, that the
 set of $e$ (resp. $q$) satisfying the arithmetical conditions
 contains about $d(n-d)$ (resp. $(n+1-d)/2$) many distinct arithmetic progressions 
 with modulus $ d(n-2)$ (resp. $d(n+1-d)$).

% set of accessible degrees $e$ and its subset of point-minimal degrees
% contain a certain number of ray congruence classes $\mod d(n-2)$
%  and similarly for the corresponding set of point numbers $q\mod d$; 
%for $e$  this implies that  the density of the accessible set  is  at least about $(n-d)/n$
%and the set of $q$ has density at least $1/2d$.
\end{itemize}
%\item
% In particular, a very general Fano hypersurface
%in $\P^n, n\geq 4$, over an uncountable base field is separably infinitely 
%rationally connected (Theorem \ref{very-general-thm}).

%*******
%In both results, $q$ is maximal for the given $e$.
%In the case $d=n$,  $q(e)=e-n+1$ and 
%an evident dimension count shows that for
%$q(e-1)+1<k\leq q(e)+1$, hence $e=n-2+k$ is indeed the minimal  
%(projective) degree of
%a rational curve through $k\geq 1$ general points. so $X$ is minimally
%rationally $k$-connected for all $k$ 
%For $d<n$ it is still true is still true for those
%drgrees $e$ that occur,  that $e$ is the minimal degree
%of a rational curve through   $k=q(e)+1$ general points so $X$
%is minimally rationally $q(e)$-conncted.
%*********** \par
To my knowledge these are the first examples (in any
characteristic) of high-degree 
balanced rational curves
on Fano hypersurfaces  
except hypersurfaces 
of very low degree or dimension $\leq 3$; ditto for the minimal rational connectivity results
(separable or not). While our results are essentially optimal for
the case $d=n$, they need not be optimal for $d<n$, indeed we see no
obstruction to the existence of balanced rational curves of all
sufficiently large degrees $e$ on a very general $X$ of degree $d<n$. \par
The proof for $d=n$, presented in \S \ref{fano},
 is based on degenerating the hypersurface to a reducible variety
$X_1\cup X_2$  called
a fan hypersurface, where $X_1$ is  a hypersurface
of degree $n$ with a point of multiplicity $n-1$ blown up at that point,
and where $X_2$ is a hypersurface of degree $n-1$. 
Using a bundle smoothing result proven in \S \ref{caudate} plus in some cases
 a vanishing theorem of Rathmann \cite{rathmann-eh-conj}, we show that
a suitable rational curve on $X_1$, glued to a union of some lines on $X_2$, is well behaved
and smooths out to a balanced rational curve on $X$.\par
The proof for $d<n$ , presented in \S \ref{d<n}, 
is based on a generalization of fan called fang where
$X_1$ and $X_2$ are blowups of a degree-$d$ hypersurface
containing, respectively, a $(d-1)$-fold $\P^{n-m-1}$ or a $\P^m$, $m\geq 2$.
\par
In \S \ref{pn} we illustrate the fan method by computing the normal bundle of a general 
rational curve in $\P^n$. \par
We begin in \S \ref{caudate} with a general result on smoothing of bundles on a curve
consisting of a 'body' together with rational tails. The general thesis is that if
the bundle is balanced on each tail and glued to the body in a sufficiently general manner,
then a smoothing is no worse, and usually better, than the bundle on the body.
Note that  a bundle on a rational tree-
including the kind envisaged in this result and its applications- need not split
as a direct sum of line bundles (see \cite{ran-normal} , Example 5.6
or Example \ref{comb}  below). Thus the  proof is not just a matter
of semi-continuity, but is rather based on bundle modifications on surfaces.\par
The two 'preliminaries' sections \S \ref{prelim1} and \S \ref{prelim2} establish
some notation, terminology and largely standard general results used in the paper.\par
An appendix by M.C. Chang analyzes the numerology of curve and hypersurface degrees arising
out of the slope-matching condition in Lemma 25.\par
I am grateful to J\"urgen Rathmann and to the referees for many helpful 
and detailed comments and references. I am greatly indebted to M. C. Chang for
her work on the numerology which constitutes the Appendix..
\section{Preliminaries}\label{prelim1}
We work over an arbitrary algebraically closed field
and use Grothendieck's projective bundle convention.
In this section we collect some definitions and facts, mostly well known,
for later reference. See \S \ref{prelim2} for more such.
\subsection{Some numerology}\label{numerology-sec}
It is convenient to set, for fixed $d\leq n$, 
\eqspl{qmax}{q_{\max}(e)=\lfloor\frac{e(n+1-d)-2}{n-2}\rfloor+1}
\eqspl{emin}{e_{\min}(q)=\lceil\frac{(q-1)(n-2)+2}{n+1-d}\rceil.}
Some trivial remarks:
\begin{itemize}\item
	For $d>3$ $e$ is point-minimal iff the remainder of $e(n+1-d)-2$ mod $n-2$ is 
	$< {n+1-d}$\item
	$q_{\max}(e)-1$ is the round-down of the slope of the normal bundle of a rational 
	curve $C$ of degree $e$  on a hypersurface $X$ of degree $d$ in $\P^n$, hence 
	$q_{\max}(e)$ equals
	the expected (i.e.  maximum possible)   number of general points
	that may lie on $C$ ,
	as well as the actual number of general points on a general $C$ if $C$ is balanced
	(note that a balanced curve is automatically unobstructed, hence lies on a unique
	smooth component of the Hilbert scheme of curves in $X$).
	\item Similarly $e_{\min}(q)$ is the expected (i.e. minimal possible)  degree
	of a rational curve that may go through $q$ general points of $X$
	and the expected degree equals the actual degree iff $X$ is minimally $(q-1)$-connected.
	\item The two functions above are mutual 'sub' -inverses:
	If $d> 3$ then $\frac{n+1-d}{n-2}<1$ so by elementary arithmetic
	$q_{\max}(e_{\min}(q))=q$
	while $e_{\min}(q_{\max}(e))\leq e$, with equality iff $q_{\max}(e-1)<q_{\max}(e)$
	iff $e$ is of the form $e_{\max}(q)$ for some $q$ iff $e$ is point-minimal.
	\item Thus, if $X$ contains a balanced rational curve of degree $e$ which
	is point-minimal,
	then $X$ is minimally rationally $(q_{\max}(e)-1)$-connected 
	(if $e$ is not point-minimal there may exist rational curves of degree $<e$
	through $q_{\max}(e)$ general points).\end{itemize}

\subsection{Fans}\label{fans}[cf. \cite{hypersurf}, sec. 4]
A 2-fan is a variety of the form $P_1\cup P_2$ where $P_1$ is a blowup
$B_p\P^n$ with exceptional divisor $E\simeq\P^{n-1}$ and $P_2=\P^n$,
so that $P_1\cap P_2$ is embedded as $E\subset P_1$ and as
a hyperplane in $P_2$. For every $d>e>0$ there is a very ample 
divisor on $P_1\cup P_2$ which is $dH-eE$ on $P_1$ and $eH$ on
$P_2$, $H=$hyperplane. A divisor of this class is said to
be of type $(d, e)$.\par
A 2-fan is the special fibre $\pi\inv(0)$ in a relative 2-fan $\pi:\cP(2)\to \A^1$
which is just $B_{(p,0)}\P^n\times \A^1$, where $P_2$ is the exceptional divisor.
The divisor $p_1^*(dH)-eP_2$ induces a divisor of type $(d,e)$ on the special
fibre and $dH$ on other fibres.\par
We will also use a generalization of 2-fans called fangs, obtained in the above 
relative setting
by blowing up  $\P^r\times 0\subset \P^n\times\A^1$ for any 
$r\leq n-2$ rather than just $(p, 0)$.
Fangs will be used implicitly  in \S 1 below to give a proof of the balancedness of a 
general rational curve of degree $e\geq n$ in $\P^n$,
and again more formally  in \S \ref{prelim2} starting in \S \ref{fangs}.
\subsection{Subvarieties of fan(g)s}
This actually applies to any ambient variety $P$ with normal crossing double points 
and no other singularities.
Let 
\[P=P_1\cup_Q P_2\]
be a transverse union of two smooth  $n$-dimensional varieties  
meeting in a common smooth
divisor $Q$ (e.g. $P$ could be a fan or fang, see \S \ref{fangs}). Let
\[X=X_1\cup_{Q_X}X_2 \subset P\]
be a subvariety consisting of smooth, codimension-$c$ subvarieties $X_i\subset P_i$
meeting $Q$ transversely in a smooth common divisor $Q_X=X_1\cap X_2$.
We have a Mayer-Vietoris sequence
\[\exseq{\O_P}{\O_{P_1}\oplus\O_{P_2}}{\O_Q},\]
and likewise for $X$.
Locally, $Q_X\subset Q$ is defined by equations $u_1,...,u_c$ on $Q$ which extend to defining
equations $u_i^j, i=1,...,c, j=1,2$ for $X_j$ on $P_j, j=1,2$.
Here the $u_i$ are part of a local coordinate system on $Q$ and
$u_i^j$ their extension to $P_j$.  Then $u_i^1$ and $u_i^2$ glue together
to a function $v_i$ on $P$ and $v_1,...,v_c$ constitute local defining equations
for $X\subset P$. In particular, $X\to P$ is a lci embedding and
\[\I_{X/P}\otimes\O_{P_j}=\I_{X_j/P_j}.\]
Consequently we have for the conormal bundles
\[\check N_{X/P}\otimes\O_{X_j}=\check N_{X_j/P_j}, j=1,2,\]
and
\[\check N_{X/P}\otimes\O_{Q_X}=\check N_{Q_X/Q}.\]
Therefore we have likewise for the normal bundles
\[ N_{X/P}\otimes\O_{X_j}= N_{X_j/P_j}, j=1,2,\]
and
\[ N_{X/P}\otimes\O_{Q_X}= N_{Q_X/Q}.\]
The latter equality
may be written suggestively as
\[N_{X/P}=N_{X_1/P_1}\cup_{N_{Q_X/Q}} N_{X_2\cup P_2}.\]
More generally, given coherent sheaves $E_i$ on $P_i$, $i=1,2$, together with
an isomorphism 
\[E_1\otimes\O_Q\to E_2\otimes\O_Q\]
there is a uniquely determined coherent sheaf $E$ on $P$ restricting to $E_i$
on each $P_i, i=1,2$. 
\subsection{Subvarieties of subvarieties of fan(g)s}\label{subsub-sec}
Notations as above, given a subvariety
\[C=C_1\cup C_2\subset X\]
where $C_1\subset X_1, C_2\subset X_2$ are smooth codimension-$c'$ 
subvarieties meeting $Q_X$  transversely in the same subvariety $Q_C$,
then $C$ is lci in $X$ with normal bundle
\[N_{C/X}=N_{C_1/X_1}\cup N_{C_2/X_2}\]
and we have an exact sequence
\[\exseq{N_{C/X}}{N_{C/P}}{N_{X/P}\otimes\O_C}.\]
Now suppose we have a smoothing, i.e. a flat family $\pi:\mathcal X\to T$ with
$\pi\inv(0)=X$ and smooth general fibre. Then there is an exact sequence
\[\exseq{N_{C/X}}{N_{C/\mathcal X}}{N_{X/\mathcal X}|_C}\]
with  $N_{X/\mathcal X}$ a trivial bundle. If $H^1(N_{C/X})=0$ and $H^1(\O_C)=0$
(e.g. $C$ is a rational tree), then $H^1(N_{C/\mathcal X})=0$ and 
$H^0(N_{C/\mathcal X})\to H^0(N_{X/\mathcal X}|_C)$ is surjective.
This means that $C$ is unobstructed as subvariety of either $X$ or $\mathcal X$
and deforms along with the smoothing $\mathcal X$. 
Now as $ C\subset X$ is an inclusion of normal-crossing
double point varieties, it induces an isomorphism $T^1_X\otimes\O_C\simeq T^1_C$.
Hence
clearly the general fibre
of the deformation of $C$ is smooth because this is true for $X$. 
Moreover given smoothings
\[\mathcal C/T\subset\mathcal X/T\subset \mathcal P/T,\]
then, possibly shrinking $T$, equations for $X\subset P$ extend to
equations for $\mathcal X\subset \mathcal P$ and likewise for
$\mathcal C$. Thus,
$\mathcal C$ is lci in $\mathcal X$ and its equations restrict to to equations of $C$ on
$X$, therefore
\[N_{\mathcal C/\mathcal X}\otimes\O_C=N_{C/X}.\]
Note that a necessary condition for a smoothing $\mathcal P$ to exist is
$N_{Q/\P_1}\otimes N_{Q/\P_2}\sim_{\mathrm {num}}\O_Q$. This is just because
\[N_{Q/P_1}=\O_{\mathcal P}(P_2)\otimes\O_Q, N_{Q/P_2}=\O_{\mathcal P}(P_1)\otimes\O_Q,\]
and $\O_{\mathcal P}(P_1+P_2)\otimes\O_Q\sim_{\mathrm{num}}\O_{Q}$. If e.g. $T=\A^1$, numerical equivalence
may be replaced by linear equivalence.

\subsection{Reducible normal bundles}\label{reducible-sec}
Normal bundles to reducible varieties in a smooth ambient space behave differently
to the case of a reducible ambient space: 
\begin{lem}
	Let $P$ be a smooth $n$-dimensional variety and $X_1, X_2\subset P$ smooth codimension-$c$
subvarieties such that $Y=X_1\cap X_2$ is a smooth divisor in each. Let $X=X_1\cup X_2$. Then
$X$ is lci in $P$ and there are exact sequences
\[\exseq{\check N_{X/P}\otimes\O_{X_1}}{\check N_{X_1/P}}{\check N_{Y/X_2}},\]
\[\exseq{N_{X_1/P}}{N_{X/P}\otimes\O_{X_1}}{N_{Y/X_1}\otimes N_{Y/X_2}}.\]
\end{lem}
\begin{proof}
Locally, $X$ has defining equations of the form
\[u_1,...u_{c-1}, u_c^{(1)}u_c^{(2)}\]
where $u_1,...,u_{c-1}, u_c^{(1)}$ are defining equations for $X_1$ and $u_c^{(2)}$
 is an equation for $Y$ on $X_1$; similarly
 $u_1,...,u_{c-1}, u_c^{(2)}$ are defining equations for $X_2$ and $u_c^{(1)}$
  This shows exactness of the
 first sequence, and the second is just its ext-dual.
	\end{proof}
Thus, unlike in the case $P$ reducible, 
the conormal $\check N_{X/P}$ (resp. normal $N_{X/P}$) restricts on $X_i$ to a
rank-1 elementary down (resp. up) modification of $\check N_{X_i/P}$ (resp. $N_{X_i/P}$)
(see \S \ref{modifications-sec}).
\begin{example}
	Suppose $X_1, X_2$ are smooth curves meeting transversely in a point $Y$.
	Then we have an exact sequence
	\[\exseq{N^0_{X/P}}{ N_{X/P}} {T^1_X}\]
where $T^1_X$, which is the usual $T^1$ sheaf of $X$,
 is a 1-dimensional vector space skyscraper  at $Y$
 and a 1st-order deformation, i.e. an element of $H^0(N_{X/P})$, is
 smoothing iff is has a nonzero image in $H^0(T^1_X)$.
 The kernel $N^0_{X/P}$ is the 'locally trivial' normal sheaf, which parametrizes
local motions of the triple $(X_1, X_2, Y)$ and fits in an exact sequence
\[\exseq{N^0_{X/P}}{N_{X_1/P}\oplus N_{X_2/P}}{Z^0}\]
where $Z^0$ is a skyscraper at $Y$ equal to
 the Zariski normal space $T_{P. Y}/(T_{X_1, Y}\oplus T_{X_2, Y})$; 
 e.g.  if $\dim(P)=2$ then $Z^0=0$. Note $N^0_{X/P}$ is not locally free on
 $X$. Anyhow the condition  $H^1(N^0_{X/P})=H^1(N_{X/P})=0 $ implies
 that $X$ has unobstructed deformations in $P$ and is smoothable in $P$ (of course the first
 vanishing implies the second for $X$ a curve).
 
  %but is the image of a locally free sheaf on $X_1\coprod X_2$.

%/*************
%which yields
%	\[0\to \carets{v_i}\to N_{X_i/P}|_Y\stackrel{j_i}{\to} N_{X/P}|_Y\to T^1_X\to 0 \]
%	where $v_1, v_2$ are tangent vectors to $X_2, X_1$ respectively at $Y$
%	and $q_1, q_2$ are 1-dimensional vector spaces.
%	Thus the images of $j_1, j_2$ are the same and correspond to locally trivial
%	local deformations of $X$, i.e. the kernel of the natural map
%	$N_{X/P}|_Y\to T^1_X$. 
%	******/
	\end{example}
\subsection{Balanced bundles}\label{balanced}\footnote{Notation here 
differs slightly from the introduction} See \cite{hypersurf} for more details.
A balanced bundle $E$ or rank $r$ on $\P^1$ has the form
\[E=r^+\O(a^+)\oplus (r-r^+)\O(a^+-1), r^+>0,\]
where the uniquely determined subbundle $r^+\O(a^+)$ is called the upper subbundle
and its rank and slope are called the upper rank and degree, respectively.
$E$ is said to be \emph{perfectly balanced} if $r^+=r$ or equivalently, $E$ is a twist of a trivial bundle.
Note that the smaller line bundle degree appearing above (i.e. $a^+-1$ if 
$r^+<r$, otherwise $a^+$) equals the round-down of the slope,
i.e. $\lfloor\deg(E)/r\rfloor$, and that $H^1(E(-t)=0$ iff $t\leq\lfloor\deg(E)/r\rfloor+1$.
The fibre of the upper subbundle at a point $p$,
which is a subspace of the fibre $E_p=E\otimes k(p)$, is called the
upper subspace at $p$.\par
A lci rational curve $C$ on a smooth variety $X$ is said to be balanced (resp. perfect)
if its normal bundle $N_{C/X}$ is balanced (resp. perfectly balanced). \par
Balancedness of $E$ is equivalent to rigidity, i.e. vanishing of
$H^1(\check E\otimes E)$, and in particular it is an open property.\par
For our purposes it is important to consider 'balanced'
bundles on rational trees. 
We could define a bundle $E$ on a rational \emph{chain} $T$ to be balanced if its
restriction on any connected subtree is a direct sum
of line bundles of total degree $a^+$ or $a^+-1$.  
This condition is satisfied whenever the restriction of $E$
on any component of the chain is balanced and the gluing at the nodes is 
general (\cite{hypersurf}, Lemma 2). 
In particular, if $T$ is a 2-component chain $T_1\cup_pT_2$ then $E$
is balanced provided $E_{T_1}, E_{T_2}$ are balanced and the upper subspaces
if $E_{T_1}$ and $E_{T_2}$ are transverse.
If $E$ is balanced on $T$ then
in any smoothing of $(T,E)$, the general fibre is balanced
(see \cite{hypersurf}, Lemma 5). A different approach to smoothing
is given in \S \ref{caudate} below.
\subsection{Modifications} \label{modifications-sec}
Given a vector bundle $E$ on a variety
$X$, a Cartier divisor $D$ on $X$, and an exact sequence
of locally-free $\O_D$-modules, where $E_D$ denotes $E\otimes\O_D$,
\[\exseq{P}{E_D}{Q},\]
the \emph{ elementary down modification} of $E$ corresponding
to $Q$ is an exact sequence
\[\exseq{M_Q(E)}{E}{Q}.\] Then $M_Q(E)$ is a locally free $\O_D$-module 
containing $E(-D)$ and
fits in another exact sequence
\[\exseq{E(-D)}{M_Q(E)}{P}\] which yields
\[\exseq{E}{M_Q(E)(D)}{P\otimes\O_D(D)}\]
Set $M^P(E):=M_Q(E)(D)$, called the \emph{elementary up modification} of
$E$ corresponding to $P$. 

Locally, if $t$ is an equation for $D$,
then there is a local basis $x_1,...,x_r$ of $E$ such that
$x_1,...,x_s$ (resp. $x_{s+1},...,x_r$) yields a local basis for
$P$ (resp. $Q$), then 
$tx_{s+1}, ...,tx_r, x_1,...,x_{s}$ is a local basis of $M_Q(E)$
(resp. $x_{s+1},..., x_r, x_1/t,..., x_s/t$ is a local basis for $M^P(E)$),
adapted in each case to the appropriate exact sequence above.
Note in our applications $Q$ will have constant rank but this is not required.\par
%For the restrictions on $D$, we also have exact
%\[0\to P\to M_Q(E)_D\to E_D\to Q\to 0.\]
For restriction on $D$, we have a exact sequences
\[\exseq{Q\otimes\O_D(-D)}{M_Q(E)\otimes\O_D}{P},\]
\[\exseq{Q}{M^P(E)\otimes\O_D}{P\otimes\O_D(D)}.\]
Thus an elementary down (resp. up) modification turns  a sub to a quotient 
(resp. a quotient to a sub).
%Set $M^P(E):=M_Q(E)(D)$, called the \emph{elementary up modification} of
%$E$ corresponding to $P$. 
\par
A \emph{modification} of $E$ is the composition of a sequence of elementary down and
elementary up modifications.
These constructions apply in particular to the case of a bundle $E$ on a
curve $C$, in which case a modification may be realized as the composition
of a single down and a single up modification (or vice versa).
For an elementary modification, corresponding to a smoothly supported reduced divisor
$D=\sum p_\l$ on $C$,  $P$ and $Q$
are just a sub and quotient vector space of $E_D=\bigoplus E\otimes k(p_\l)$.
The modification is said to be \emph{general} if $D$ is reduced and the sub
or quotient in question are general.	
If $D$ is supported on a unique component $F$ of $C$
and $E$ restricted on  $F$ is balanced,
the modification is said to be \emph{in general position} 
if the induced map from the upper subbundle (see \S \ref{balanced})
\[E_F^+\otimes \O_D\to Q\] has maximal rank .
\begin{lem}\label{mod-lem}
	Assumptions as above, let $E$ be a bundle on a curve $C$ and
	$F$ a component of $C$ such that
 $E_F$ is balanced  with upper rank $r^+$ and upper degree $a^+$ 
	and let $E'=M_Q(E)\subset E$ be 
	an elementary modification of colength $s$ in general position
	cosupported on $F$. Then
	if $s<r^+$, we have 
	\[r^+(E')=r^+-s, a^+(E')=a^+(E).\] Otherwise,
	\[r^+(E')=r+r^+-s, a^+(E')=a^+-1.\]
\end{lem}
\begin{proof}
	This follows easily from the fact that the induced map
	$E_F^+(p)\to sk(p)$ has maximal rank by generality (here and elsewhere,
	$sA$ where $s\in\N$, denotes $\bigoplus\limits_1^sA$ wherever this makes sense).
\end{proof}
%A similar argument proves the following
%\begin{lem}\label{torsion}
%	Let $E$ be a balanced bundle and $\tau$ a torsion sheaf on a smooth curve $C$, and let
%	$\phi: E\to\tau$ be a sufficiently general surjection. Then $\ker(\phi)$ is balanced.
%	\end{lem}
%	\begin{proof}
%		Using a Jordan-Holder filtration we may assume $\tau$ has length 1, i.e.
%		$\tau=k(p), p\in C$. Then the assertion follows from Lemma \ref{mod-lem}.
%		\end{proof}

%/**************
%For a modification in general position it is easy to check,
%setting $s=\rk(Q)$,
%that $M_Q(E)|_F$ is also balanced, with upper rank $r^+-s$ if $r^+>s$
%or $r+r^+-s$ if $r^+\leq s$. See Lemma \ref{mod-lem}
%for a similar statement for point modifications and Lemma \ref{mod} below for a more general result.\par
%Unless otherwise stated, all elementary modifications we use in the curve case will be 
%of the pointwise variety.
%A pointwise elementary down (resp. up) modification of colength $s$
% of a bundle  $E$ on a  curve $C$  is just an inclusion $E'\subset E$
% (resp.  $E\subset E"$)
%such that $E/E'\simeq sk(p)$ (resp. $E"/E\simeq sk(p)$) for some smooth point $p\in C$.
%********** / \par

\subsection{Blowing up normal bundles}
Elementary modifications occur often in the geometry of embedded curves.
One example is the following standard result which to save notation
we have stated just for a curve $C$ but with evident modifications is equally valid
for $C$ any lci subvariety (which will naturally get blown up in the blowup 
of $X$).
\begin{lem}\label{blowup-lem}
	Let $C$ be a lci curve on a smooth variety $X$ and let $Y$ be a lci 
	subvariety of codimension $s$ in $X$ meeting $C$ schematically in a 
	unique point $p$ smooth on $C$.
	Let $\pi:X'\to X$ be the blowup of $Y$ with exceptional divisor $E$ 
	and let $C'$ the birational transform of $C$
	on $X$. Then \par
	(i) $N_{C'/X'}$ is the elementary down modification
	of colength $s-1$
	of $N_{C/X}$ corresponding to the image of $T_pY$ in $N_{C/X}(p)$.\par
	(ii) Under the natural identification of $(N_{C'/X'})_p$ with $T_pE$, 
	the kernel of the map $(N_{C'/X'})_p\to (N_{C/X})_p $ coincides with the vertical subspace of
	$T_pE$, i.e. the tangent space to the fibre of $E\to Y$.
	\end{lem}  
\begin{proof}For convenience we work with conormal bundles denoted $\check N_{C/X}$ etc.
	If $\pi:X'\to X$ denotes the blowup map, $\pi^*\check N_{C/X}$ is clearly
	a subsheaf of $\check N_{C'/X'}$ and coincides with it locally off $p$,
	so it suffices to identify the image at $p$.
	We can choose local coordinates at $p$ of the form
	$y, x_1, ...,x_{s-1}, x_{s},..., x_n$ so that 
	$y$ defines $p$ on $C$, $x_1,...,x_n$ define $C$ and
	$y, x_1,...,x_{s-1}$ define $Y$. Then $x_1,..., x_n$ yield a local basis for
	$\check N_{C/X}$ while, in a suitable affine open in $X'$
	containing $C'$ with coordinates $y, x_1/y, ...,x_{s-1}/y, x_s,...,x_n$,
	a basis for $\check N_{C'/X'}$ is $x_1/y,...,x_{s-1}/y, x_s,...,x_n$.
	This proves the dual statement for conormals which is 
	equivalent to assertion (i).\par
	As for (ii), it follows from the above computation or, 
	just as  well, from the diagram
	\[
	\begin{matrix}
	T_pE&\simeq&(N_{C'/X'})_p\\
	\downarrow&&\downarrow\\
	T_pY&\subseteq& (N_{C/X})_p.	
	\end{matrix}
	\]
	
	\end{proof}
\section{Bundles on caudate curves}\label{caudate}
The purpose of this section is to prove a general and elementary result about smoothing
of vector bundles on curves endowed with multiple tails. This result
permits construction of some balanced vector bundles on rational curves and in particular
to prove the existence of some balanced rational curves.
The result is stated in much greater generality than is needed
for the applications to minimal rational curves given in this paper, in the hope that
it might enable further such applications.
See also \cite{coskun-scrolls}, \cite{smith-rat-trees}, \cite{ran-normal}, \cite{hypersurf}  for other results
on bundles on rational trees.
\par
By definition, a \emph{rational tree} is a nodal curve that is a tree of smooth
rational curves. A \emph{broken comb} is a connected nodal curve of the form
\[ C=B\cup \bigcup T_i\]
where $B$, the \emph{base} (aka the body), is a connected 
nodal curve and each \emph{tooth} (aka tail) $T_i$ is a rational tree 
meeting $B$ in a unique smooth
point called its \emph{root} and meeting no other $T_j$. A broken comb is \emph{rational}
if $B$ is a rational tree. A \emph{rational comb} is a rational broken comb that is unbroken,
i.e. where $B$ and each $T_i$ are $\simeq\P^1$.\par
Unlike the irreducible case, or for that matter the case of
rational \emph{chains}- see \cite{hypersurf}- 
even nice bundles on rational combs need not split
as direct sums of line bundles. The following example is essentially taken from
\cite{ran-normal}.
\begin{example}\label{comb}
	Let $C=B\cup\bigcup\limits_{i=1}^t T_i$ 
	be a  rational comb and let $E$ be a vector bundle on $C$
	whose restriction on each $T_i$ is isomorphic to $\O\oplus \O(-1)$,
	with $E_B$ an arbitrary vector bundle and with  general gluing at nodes.
	Then $h^0(\check E\otimes E)\geq t\ $ hence, if $t\geq 5$, then
	$h^0(\check E\otimes E)>4=\chi(\check E\otimes E)$, therefore 
	$h^1(\check E\otimes E)>0$. Consequently, $E$ is not  balanced
	and in fact not a direct sum of line bundles. Nonetheless, as we shall see below in Example \ref{comb2},
	 Theorem \ref{main-thm} below applies to $E$, showing that a smoothing of $E$ is a deformation
	of a general down modification of $E_B$ at the nodes, hence is as well-behaved as possible.
%	Consequently, if $E_B\simeq\O(a_1)\oplus\O(a_2)$ then $E'\simeq \O(b_1)\oplus\O(b_2)$
%	with $|b_1-b_2|\leq\max(|a_1-a_2|-t, 1)$. Informally, attaching an $\O\oplus\O(-1)$ tail
%	works like an elementary down modification.
%	
	
\end{example}
%**************
%A bundle $E$ on a rational tree $T$ is said to be \emph{strongly balanced}
%if there exists a subbundle $E^+\subset E$ whose restriction on each component
%$U$ of $T$ coincides with the upper subbundle of the restriction $E_U$, i.e.
%\[(E^+)_B=(E_B)^+.\]
%
%
%If $B$ is a nodal curve, a \emph{blowup} of $B$ is, by definition if you will, a 
%nodal curve $B'$ equipped with a surjective map 
%$g:B'\to B$  whose non-point fibres are rational trees, each mapping
%to some node of $B$ (the set of nodes involved is called the \emph{support} of the blowup). 
%In this case, if $E_B$ is a vector bundle
%on $B$, we will use $E_{B'}$ to denote $g^*E_B$. Note that any smoothing
%of $B$ is dominated by a smoothing of $B'$ with smooth total space,
%where $B'$ is some blowup of $B$.
%*************\par
 If $E_C$ is a bundle on a curve $C$, by a
\emph{smoothing} of $(C, E_C)$ is meant an irreducible surface $S$ endowed with a flat map to a smooth curve 
$S\to T$ with smooth general fibre and special fibre $C$ (with multiplicty 1), 
plus a vector bundle $E$ on $S$
that restricts to $E_C$ on $C$. Similarly, for a  subset $A\subset C$, a \emph{partial
	smoothing} at $A$ assumes only that there is a neighborhood $U$ of $A$ on $S$
such that the intersection of the general fibre with $U$ is smooth.
\begin{thm}\label{main-thm}
	Let $C=B\cup\bigcup T_i$ be a broken comb with teeth 
	$T_1,...,T_k$ and respective roots $p_1,...,p_k$
	and let $E_C$ be a vector bundle on $C$.
	Assume\par (i) on each component of each $T_i$, $E_C$ is balanced;\par
	and for each $T_i$ either\par
	(ii) the gluing at each node $q$ of $C$ on $T_i$ of  
	the restrictions of $E_C$ on the components of $C$ through $q$  is 
	sufficiently general; or\par
	(ii)' the restriction of $E_C$ on  $T_i$ is perfectly balanced, i.e. 
	 a twist of a trivial bundle,
	i.e. has the form $\mathrm{(vector\  space)}\otimes\mathrm{(line\  bundle)}$.\par
	Then any partial smoothing of $(C, E_C)$ at $\bigcup T_i$ is the pullback by a 
	birational map of a deformation
	of a general modification of some twist $E_B\otimes\O_B(\sum m_ip_i)$ at $p_1,...,p_k$.
%	\par	
%	Moreover the following numerical  relation holds:
%	\eqspl{degree-eq}{
%		c_1(E_{T_i})=rm_i+r_i, i=1,...,k	
%	}
%	where $r=\rk(E_C)$ and, with the above notation, $r_i=\rk(P)_{p_i}$ for an up modification and
%	$r_i=-\rk(Q)_{p_i}=-r+\rk(P)_{p_i}$ for a down modification. Also, $P$
%	coincides at $p_i$ with the upper subspace of $E_{T^0_i}$ where $T^0_i$ is the unique component 
%	of $T_i$ through $p_i$.
\end{thm}
\begin{rem}
	We recall that a general (down) modification of a bundle $E$ at a smooth point $p$
	on a curve means the kernel of a map $E\to Q$ to general quotient $Q$ of $E\otimes\kk(p)$.
	Ditto for up modification and several points.
	\end{rem}
\begin{rem}
	Here the genus of $C$ is arbitrary but in applications it will be zero.
	In fact, the only case used here is that of a rational comb, i.e. where $B$ and each $T_i$
	are $\P^1$.\par
	\end{rem}
\begin{rem}
 Note that the Theorem applies an \emph{arbitrary} 1-parameter
	partial smoothing rather than just
	a 'sufficiently general' one or, for that matter, a multi-parameter smoothing
	dominating a versal deformation of the curve, where the nodes smooth independently. 
	This feature  is crucial for applications to curves on fans
	because when the curve smooths together with the fan, the nodes lying on the
	fan's double locus smooth \emph{simultaneously}, so this smoothing of the curve is 
	never general.
%	In case $r_i=r$ the elementary modification at $p_i$ is just a twist can may
%	my lumped into the twisting divisor $\sum m_ip_i$.
\end{rem}
\begin{rem}
	Below we will use the fact that in a partial smoothing, any extremal fibre
	component contained in the smooth part of the total space must be a (-1) curve.
	This follows easily from the fact that the fibre is reduced (see below). 
	\end{rem}
\begin{rem}
	Regarding the meaning of 'general gluing' in terms of the inductive procedure
	used  in the proof below.
	This proceeds 'from the outside in' (i.e. toward the 'spine' $B$).
	In practice, general gluing at a node $q$ lying on an 'outer' component $F$
	and an 'inner' one $F^*$ means that
	the upper subspace of $E_{F^*}$ at  is transverse to the 'distinguished subspace' of 
	$E_F$ at $p$, the latter being the upper subspace  of a modification of the original $E_F$
	previously constructed in the course of the proof, which has to do with upper subbundles
	on components 'further out' than $F$. It is possible that $E_{F^*}$ is a twist 
	of a trivial bundle, so its upper subspace is all of $E_p$, but after the next modification
	of $E$ the resulting bundle will generally be nontrivial on $F^*$. 
\end{rem}

\begin{proof}[Proof of theorem]
	Given a partial smoothing $(E, S)$, we first resolve all singularities of the surface
	$S$ lying on
	$\bigcup T_i$ to obtain a smoothing with  total space that is smooth in a neighborhood
	of the preimage of $\bigcup T_i$, this at the cost
	of  augmenting the $T_i$ by some further rational trees $K_j$
	on which $E$
	is trivial.  Suitably refreshing notation, we have a
	fibred surface $\pi:S\to B$ with \[\pi\inv(0)=B\cup\bigcup T_i\cup\bigcup K_j\]
	such that $S$ is smooth along $\bigcup T_i\cup\bigcup K_j$, 
	and such that each restriction  $E_{T_i}$  balanced
	and each $E_{K_j}$ is trivial. Moreover, each ${K_j}$ meets
	$B\cup\bigcup T_i$ in exactly 2 points $p'_j, p"_j$ at most one of which is on $B$.
	We may write $E_{K_j}=U\otimes\O_{K_j}$ for a vector space $U$. Then letting
	 $Z', Z"$ denote the components other than $K_j$
	through $p'_j, p"_j$ respectively, the generality hypothesis (ii) above
	implies that the identification of the fibres  $E_{p_j'}$ and $E_{p"_j}$ with
	$U$ is general and in particular the upper subspaces of $E_{Z'}$ at $p_j'$
	and $E_{Z"}$ at $p"_j$, considered as subspaces of $U$, are mutually
	in general position.\par
	Now the proof is an inductive procedure on the  irreducible components of 
	the 'multitail' $\mathcal T= \bigcup T_i\cup\bigcup K_j$,
	proceeding 'inward' towards $B$, where each step eliminates by contraction
	an extremal component. The procedure works separately for each $T_i$ and we will
	work on a $T_i$ for which (ii) above holds as the case of (ii)' is similar and simpler. 
	We start with the initial step, which is essentially identical
	to the inductive step. Let $F$ be an extremal component of $\mathcal T$, i.e.
	$F$ meets the rest of the curve, say $G$,  in a single point $p$. Because $G$ is reduced
	and $F.(F+G)=F.(\mathrm{fibre})=0$, we have $F^2=-F.G=-1$ so $F$ must be a $(-1)$
	curve (initially $F$ must be a component of some $T_i$ but this is unimportant). 
	By assumption we can write
	\[E_F\simeq r^+\O_F(d^+)\oplus (r-r^+)\O_F(d^+-1).\]
	Replacing $E$ by its twist $E(d^+F)$, we may assume $d^+=0$. Now
	if $r^+=r$, i.e. $E_F$ is a twist of a trivial bundle, we may as well
	assume $E_F\simeq r\O_F$.
	If $r^+<r$, perform an elementary down modification on $E$ corresponding to the quotient
	\[E\to (r-r^+)\O_F(-1).\]
	 This modification yields a subsheaf $E'\subset E$, equal to $E$
	off $\pi\inv(0)$, with $E'_F\simeq r\O_F$. Moreover if $F^*$ is the unique other component
	of $\pi\inv(0)$ through $p$ then $E'_{F^*}$ is the elementary modification of $E_{F^*}$ at $p$
	corresponding to the corresponding pointwise quotient of vector spaces
	\[E(p)\to (r-r^+)\O_F(-1)(p).\]
	Furthermore, if $G$ is the subcurve of $\pi\inv(0)$
	complementary to $F\cup F^*$ then $E'_G=E_G$.\par
	 Now if the aforementioned  
	 $F^*$ is a component of some $T_i$ then by our general gluing hypothesis,
	$E'_{F^*}$ is balanced and its upper subspace at all nodes on $F^*$, being a sum or intersection
	of 2 general subspaces (compare Lemma \ref{mod-lem}), is general. \par 
	In the other case, $F^*$ is a component
	of some $K_j$ which meets $B\cup\bigcup T_i$ in another point $p'=K_j\cap F^{**}$
	where $F^{**}$ is a component of $B\cup\bigcup T_i$
	and by general gluing plus the fact that $E_{K_j}$ is a trivial bundle $U\otimes_\C\O_{K_j}, U\simeq\C^r$,
	 ensures as above that
	$(E_F)^+(p)$ and $(E_{F^{**}})^+(p')$, as subspaces of $U$, are general and meet
	transversely. In particular, $E'_{F^*}$ is balanced.

%****************	
%which, thanks to our general gluing hypothesis, is general.
%	Hence again we may as well assume
%	$E_F\simeq r\O_F$.
%	************
	\par
	Next, we blow down $F$. Set $\O_{rF}=\O_S/\O_S(-rF)$ and consider the standard exact sequence
	\[\exseq{E'\otimes\O_F(-rF)}{E'\otimes\O_{(r+1)F}}{E'\otimes \O_{rF}}.\]
	Using the fact that $E'_F=r\O_F$ and
	 $\O_F(-F)=\O_F(1)$ it follows easily that if $\hat F$ denotes the
	formal completion of $S$ along $F$, then \[E'\otimes \O_{\hat F}\simeq r\O_{\hat F}.\]
	Consequently if we let $f:S\to S_1$ denote the blowing down of the $(-1)$ curve $F$,
	then by the formal function theorem 
	$f_*(E')$ is locally free near $q=f(F)$ (also $R^1f_*(E')=0$).
	Hence if we let
	\[E_1=f_*(E')\] then identifying the general fibre $Y$ of $S/B$ and $S_1/B$,
	we have $(E_1)_Y\simeq E'_Y$. Thus we get a modified family with similar properties whose special 
	fibre has a smaller multitail, and we may continue the process until there is no multitail,
	which is the desired family. 
%	Finally the numerical relation
%	\eqref{degree-eq} is immediate from conservation of degree.

\end{proof}
%\begin{rem}
%Regarding the meaning of 'general gluing' in terms of the inductive procedure in the proof.
%This proceeds 'from the outside in' (i.e. toward the 'spine' $B$).
%	In practice, general gluing at a node $q$ lying on an 'outer' component $F$
%	and an 'inner' one $F^*$ means that
%	the upper subspace of $E_{F^*}$ at  is transverse to the 'distinguished subspace' of 
%	$E_F$ at $p$, the latter being the upper subspace  of a modification of the original $E_F$
%	previously constructed in the course of the proof, which has to do with upper subbundles
%	on components 'further out' than $F$. It is possible that $E_{F^*}$ is a twist 
%	of a trivial bundle, so its upper subspace is all of $E_p$, but after the next modification
%	of $E$ the resulting bundle will generally be nontrivial on $F^*$. 
%\end{rem}
\begin{rem}
	Similar results in the case of a rational chain, rather than broken comb, or a comb with few teeth
	are given in \cite{hypersurf}, \S 1.
	\end{rem}
\begin{cor}\label{tree-cor}
	Let $T$ be a rational tree and let $E_T$ be a vector bundle on $T$ 
	such that for each component $S$ of $T$ either\par
	(i) the restriction $E_S$ is balanced and the gluing at each node 
	on $S$ is general; or\par
	(ii) $E_S$ is a twist of a trivial bundle.\par
	 Then any smoothing of $(T, E_T)$ has balanced general fibre. 
\end{cor}
\begin{rem}
	Case (ii) above is 'almost' a special case of case (i), except when there
	is a pair of intersecting components $S_1, S_2$ with $E_{S_1}, E_{S_2}$
	trivial; then the gluing of the trivial bundles $E_{S_1}, E_{S_2}$ at
	$S_1\cap S_2$ is immaterial. 
	\end{rem}
Note that by Example \ref{comb2} below,
it is possible under the hypotheses of the Corollary  to have
$h^1(\check E_T\otimes E_T)>0$, a condition which for $\P^1$ is equivalent
to non-balancedness.
The Corollary may be used in lieu of Lemma 2 or Lemma 7 of \cite{hypersurf}
to show existence of some  balanced rational curves of low degree $e$ 
on general Fano hypersurfaces
of degree $d\leq n$, and will be used for a similar purpose in \S 3 below
for the case $d=n$ and $e\geq n-1$.
\begin{cor}
	Let $f:X\to S$ be a proper flat family of nodal-or-smooth curves
	with general fibre isomorphic to  $\P^1$, over an irreducible variety $S$. 
	Let $\partial S\subset S$ be the locus of singular fibres.
	Let $E$ be a vector bundle on $X$. Suppose that $T:=f\inv(s_0)$ together with
	$E_T$ satisfy the hypotheses of Corollary \ref{tree-cor}. Then 
	there is a neighborhood $U$ of $s_0$ in $S$ such that for every $s\in U\cap(S\setminus\partial S)$,
	%$s_0$ is not in 
	%the closure of the locus of  $s\in S$ such that $f\inv(s)\simeq\P^1$
	%and $E_{f\inv(s)}$ is not balanced. 
	%	\[s_0\not\in\supp(R^1f_*(\check E\otimes E)).\]
	$E_{f\inv(s)}$ is balanced; equivalently,
	\[\mathrm{supp}(R^1f_*(\check E\otimes E))\cap U\subset \partial S\cap U.\]
\end{cor}
The Corollary is interesting because it applies in situations where
standard semi-continuity fails because, with the
above notation, one has $H^1(\check E_T\otimes E_T)\neq 0$- see Example \ref{comb} below. 
Then we conclude that $R^1f_*(\check E\otimes E)$ is nontrivial and 
locally supported on the boundary.
%In that case
%the Corollary implies that the support of $R^1f_*(\check E\otimes E)$
%is contained in $\partial S\cap U$.
\par
Returning to the general situation of the Theorem,  it actually
implies more, namely to the effect that, when nontrivial modification
get involved, a general smoothing
of $(C, E)$ is 'better behaved' than $E_B$. To make this precise, 
it is convenient to use the language
of partitions. Suppose $E$ is a vector bundle on $\P^1$of the form

\[ E\simeq\bigoplus\limits_{i=1}^sr_i\O(d_i),\ \  d_1>d_2>...>d_s.\] 
The subbundles
\[E_j=\sum\limits_{i=1}^j r_i\O(d_i)\]
are canonically defined and form the Harder-Narasimhan filtration of $E$:
\[E_1\subset E_2\subset....\subset E_s=E.\] 
We associate to $E$ the partition $\Pi(E)$ with blocks of height $d_i$
and width $r_i$,
$i=1,...,s$ and total width $r$, which may be viewed as usual as a subset
of the first quadrant in $\R^2$ containing the origin.
In addition to being partially ordered by inclusion,
these partitions are lexicographically ordered via 
the degree sequence $(d_i)$ and if $E'$ is a general
member of a deformation of $E$ then \[\Pi(E')\leq\Pi(E).\] Given  a partition $\Pi$
of degree $d$ and width $r$ and an integer $k$, the \emph{elementary modification}
of type $k$ of $\Pi$, denoted $M_k(\Pi)$,
is the lexicographically smallest partition $\Pi'$ of width $r$ and degree $d+k$, 
such that \[\Pi'=\Pi \ \ \mathrm{if}\ \  k=0,\]  \[  \Pi'\supset \Pi, \ \ \mathrm{if}\ \  k>0,\ \ \] and
\[ \Pi'\subset\Pi, \ \ \mathrm{if}\ \  k<0.\]
One way to define $M_k(\Pi)$  is inductively as $M_1(M_{k-1}(\Pi))$ ($k>0$)
or $M_{-1}(M_{k+1}(\Pi))$ ($k<0$),
where $M_1(\Pi)$ (resp. $M_{-1}(\Pi)$) replaces the first (resp. last)  
column of height $d_r$ (resp. $d_1$) by a 
column of height $d_r+1$ (resp. $d_1-1$). %Likewise in the case $k<0$.
\par
A  modification corresponding to  $E\to Q=\bigoplus\limits_{\l=1}^t Q_{p_\l}$ 
supported on $D=\sum p_\l$ is said to be 
\emph{in general position} if for each $i$ the induced map
\[E_{i}\otimes\O_D\to Q\] has maximal rank. The following is a 
simple generalization of Lemma \ref{mod-lem}.
\begin{lem}\label{mod} If $E'$ is an elementary modification in general position of $E$
	(up or down, at one or more points), 
	and \[\deg(E')=\deg(E)+k\] then \[\Pi(E')=M_k(\Pi(E)).\]
\end{lem}
\begin{proof}It suffices  treat the case of a down modification. Let $j$
	be smallest such that $E_j\otimes\O_D\to Q$ is surjective. Then there is
	an exact sequence
	\[\exseq{\bigoplus\limits_{i<j}r_i\O(d_i-1)\oplus r'_j\O(d_j-1)}{E'}{(r_j-r'_j)\O(d_j)
		\oplus\bigoplus\limits_{i>j}r_i\O(d_i)} \]
	with $r_1+...+r_{j-1}+r'_j=\l(Q)$, $0<r'_j\leq r_j$.
	Such a sequence automatically splits and this suffices to imply that $E'$ has the
	desired partition.
\end{proof}
Therefore the Theorem implies (compare \cite{hypersurf}, Lemma 7):
\begin{cor}
	Assumptions as in the Theorem, if 
	$B\simeq\P^1$ and
	$(C', E')$ is a smoothing of $E$ then
	$\Pi(E')\leq M_k(\Pi(E_B))$, where $k=\sum\deg(E_{T_i})$.
\end{cor}
% \begin{rem}
% 	Suitably interpreted in terms of the Harder-Narasimhan filtration on $E_B$, 
% 	the Corollary should hold for $B$ of any genus and should imply in some cases that
% 	the smoothing in stable.
% 	\end{rem}

% 
%********* \par
%The Theorem is applicable more generally to show that when $E$ is balanced
%on the tails then a smoothing will 'flatten out'
% the Harder-Narasimhan (HN) filtration of $E_B$. The general statement is notationally
%messy so we will
%just state a sample special case.
%\begin{cor}
%	Let $C=B\cup T$ be a nodal curve where $B$ is smooth and $T\simeq \P^1$.
%	Let $(d.(E_B), r.(E_B))$ be the degrees and ranks of the HN graded pieces of
%	$E_B$ in descending slope order.  % where $d_i(E_B)/r_i(E_B)>d_{i+1}(E_B)/r_{i+1}$. Assume 
%	Assume $E_T\simeq r_+\O_T\oplus (r-r_+)\O_T(-1)$
%	and that $\frac{d_1(E_B)-r_+}{r_1(E_B)}>\frac{d_2(E_B)}{r_2(E_B)}$.\par
%	Then a smoothing $E'$ of $E_C$ has $d_1(E')=d_1(E_B)-r_+, r_1(E')=r_1(E_B)$.
%	\end{cor}
%******************
\begin{example}\label{comb2}[Example \ref{comb} revisited]
%	Let $C=B\cup\bigcup\limits_{i=1}^t T_i$ 
%	be a  rational comb and let $E$ be a vector bundle on $C$
%	whose restriction on each $T_i$ is isomorphic of $\O\oplus \O(-1)$,
%	with general gluing at nodes.
%	Then $h^0(\check E\otimes E)\geq t$ hence if $t\geq 5$ then
%	$h^0(\check E\otimes E)>4=\chi(\check E\otimes E)$ hence 
%	$h^1(\check E\otimes E)>0$. Consequently, $E$ is not by any reasonable 
%	definition balanced
%	(in fact $E$ is not a direct sum of line bundles). Nonetheless,
	Notations as in Example \ref{comb}, Theorem \ref{main-thm} applies to $E$, showing that a smoothing $E'$ of $E$ is a deformation
	of a general down modification of $E_B$ at the nodes.
	Consequently, if $E_B\simeq\O(a_1)\oplus\O(a_2)$ then $E'\simeq \O(b_1)\oplus\O(b_2)$
	with $|b_1-b_2|\leq\max(|a_1-a_2|-t, 1)$. Informally, attaching an $\O\oplus\O(-1)$ tail
	works like an elementary down modification.

\end{example}

\section{Curves in projective space}\label{pn}
Here as a warmup for fan-like methods we will prove the well-known fact
(see \cite{ran-normal} for a longer proof ):
\begin{prop}\label{pn-prop}
	A general rational curve of degree $e\geq n$ in $\P^n$ is balanced.
	\end{prop}

\begin{proof}
	\ul{Case 1} (probably well known): a rational normal curve $C\subset \P^n$.\par
	 We must show $N_{C/\P^n}=(n-1)\O(n+2)$.\par
	\emph{Proof 1} (arbitrary char.): We work inductively, the cases $n=1, 2$ being easy. 
	Degenerate $C$ to $C_0=C'\cup_p L$ where $C'\subset\P^{n-1}$
	is rational normal and $L$ is a transversal line. See the discussion
	in \S \ref{reducible-sec}.
	By induction, we have $N_{C'/\P^{n-1}}=(n-2)\O(n+1)$. The natural exact sequence
	\[\exseq{N_{C'/\P^{n-1}}}{N_{C'/\P^n}}{\O(1)|_{C'}}\]
	splits because $C'$ is the intersection of a minimal scroll in $\P^n$ with $\P^{n-1}$,
	hence \[N_{C'/\P^n}=(n-2)\O(n+1)\oplus\O(n-1).\]
	Then $N_{C_0/\P^n}|_{C'}$ is the rank-1 elementary up modification
	of $N_{C'/\P^n}$ at $p$ corresponding to $T_pL$ which is not tangent to $\P^{n-1}$, hence
	  \[N_{C_0/\P^n}|_{C'}=(n-2)\O_{C'}(n+1)\oplus\O_{C'}(n),\] with upper subbundle
	  $N_{C'/\P^{n-1}}=(n-2)\O_{C'}(n+1)$.
Similarly,  
\[N_{C_0/\P^n}|_L=\O_L(2)\oplus (n-2)\O_L(1).\]
Now the natural maps
\[i_1: N_{C'/\P^n}\to N_{C_0/\P^n}\otimes\O_{C'}, i_2: N_{L/\P^n}\to N_{C_0/\P^n}\otimes\O_L\]
have the same codimension-1 image at $p$, namely the hyperplane dual to
the element of $\check N_{C_0/\P^n}(p)$ coming from 
the unique order-2 minimal local generator
of $\I_{C_0}$, i.e. a generator of
$\ker(\check N_{C_0/\P^n}\to \check N_{C'/\P^n})=\ker(\check N_{C_0/\P^n}\to \check N_{L/\P^n})$.
%
%which is also essentially the unique minimal generator of $\I_{C_0}$ that
%is not a minimal generator of $\I_{C'}$ or $\I_L$, i.e generates the kernels of
%$\I_{C_0}\otimes\kk(p)\to\I_{C'}\otimes\kk(p)$ and $\I_{C_0}\otimes\kk(p)\to\I_{L}\otimes\kk(p)$.
 Considering $i_1(p)$, the image in question coincides with the upper subspace
of $N_{C_0/\P^n}|_{C'}$ at $p$, coming from the $(n-2)\O(n+1)$ (unique) subbundle.
On the other hand, the upper subspace of $N_{C_0/\P^n}|{L}$, coming from the $\O(2)$,
is clearly not in the image of $i_2(p)$.
Since the images of $i_1(p)$ and $i_2(p)$ coincide, it follows the respective upper
	subspaces at $p=C'\cap L$  are different, i.e.  in general position, so
	$N_{C_0/\P^n}\simeq (n-1)\O_{C_0}(n+2)$. 
	\par
\emph{Proof 2} (quicker but maybe valid only in char. 0):	The normal bundle $N=N_{C/\P^n}$ has
 degree  $(n-1)(n+2)$ and rank $n-1$. 
		Textbooks (\cite{harris-ag} or\cite{griffiths-harris} p.12) show that 
		there is a unique $C$ through $n+3$ general points. Carefully examining this construction
		(or else using char. 0) shows that the locus of rational normal curves through $n+3$
		general points is a reduced singleton.		
%		Generic smoothness shows that the locus
%		of rational normal curves through 	 $n+3$ general points
%		is reduced while the textbook synthetic construction of $C$ through those points
%	
%		this locus is set-theoretically a point. 
		Consequently, this locus has trivial tangent space, i.e.
		$H^0(N(-\sum\limits_{i=1}^{n+3}p_i))=0$.
		Hence $N$ contains no line bundle of degree $n+3$ or more, so
		$N\simeq (n-1)\O(n+2)$. \par
		\ul{Case 2}:
		 $n<e<2n$.\par
		  Consider the blowup $\cX$ of $\P^n_1\times\A^1$
		in $\P^{e-n}_1\times 0$ ($\P^a_b$ is a copy of $\P^a$), 
		with natural maps $\pi:\cX\to\A^1, f:\cX\to\P^n_1$. Then
		\[X_0:=\pi\inv(0)=X_1\cup X_2,\]
		where
		\[X_1=B_{\P^{e-n}_1}\P^n_1, X_2=B_{\P^{2n-e-1}_2}\P^n_2, Z:=X_1\cap X_2= 
		\P^{e-n}_1\times\P^{2n-e-1}_2.    \]
		($\P^n_1, \P^n_2$ are copies of $\P^n$ and likewise for their subspaces.)
		$\cX$ is endowed with a relative hyperplane bundle $H=f^*\O(1)(-X_2)$
		which restricts on the general fibre of $\pi$ to $\O(1)$, on each $X_i$ to
		$\O(1)(-Z)$ and on $Z$ to $\O(1,1)$.\par
		Let $C'_1\subset X_1, C'_2\subset X_2$ be respective
		proper transforms of curves 
		$C_1$, a rational normal curve in $\P^n_1$ and
		$C_2$, a rational normal curve in its span $S\simeq\P^{e-n+1}$,
		% \bar C_{e-n+1}$ of the indicated degrees
		where $C_1$ (resp. $C_2$) meets $\P^{e-n}_1$ (resp. $\P^{2n-e-1}_2$), transversely in 1 point,
		so that $C_0=C'_1\cup C'_2\subset X_0$ is a connected nodal curve.
		We also assume $S$ is transverse to the blowup center $\P^{2n-e-1}_2$.
		Note that the 'degree', i.e. $H$-degree, of $C_1'\cup C_2'$, is $e$, e.g. because 
		$C'_1\cup C'_2$ projects in $\P^n_1$
		to the union of $C_1$, of degree $n$, and the projection $\bar C_2$ 
		of $C_2$ from $\P_2^{2n-e-1}$ which meets it in 1 point so that $\bar C_2$
		has degree $e-n$. The family of curves $C$ in $\P^n$ that we construct below degenerates to
		$C_1\cup C_2\subset X$, and that degeneration dominates a degeneration of $C$ in $\P^n$ to
		$C_1\cup \bar C_2$, so $C$ has degree $e$. Also see Example \ref{curve-degree}.

%		Now $N_{C_1/\P^n_1}\simeq (n-1)\O(n+2)$. Similarly,
We analyze $C_2$ first. Note
			\[N_{ C_{2}/S}=(e-n)\O(e-n+3)\]
		hence
		\[N_{ C_{2}/\P^n_2}\simeq (e-n)\O(e-n+3)\oplus (2n-e-1)\O(e-n+1).\]
		
%/*******	%		\[N_{ C_{2}/S}=(2n-e-2)\O(2n-e+1)\]
%			hence
%\[N_{ C_{2}/\P^n_2}\simeq (2n-e-2)\O(2n-e+1)\oplus (e-n+1)\O(2n-e-1).\]
%*************/\par

The latter bundle is not balanced, however after the blowup 
of the transverse
$\P^{2n-e-1}_2$, we get (see Lemma \ref{blowup-lem})
	\[N_{ C_{2}/\P^n_2}\simeq (e-n)\O(e-n+2)\oplus (2n-e-1)\O(e-n+1),\]
%\[N_{C_2'/X_2}\simeq (2n-e-2)\O(2n-e)\oplus (e-n+1)\O(2n-e-1)\]
which is balanced and whose upper subspace at $p$ comes from the first summand and,
identifying the normal space at $p$ with $T_pZ$, 
corresponds to the $\P_1^{e-n}$ factor, i.e.  kernel of 
the projection $T_pZ\to T_p\P_2^{2n-e-1}$ 
(see Lemma \ref{blowup-lem}) , which is also the kernel of the differential of the
blowdown map $X_2\to\P^n_2$ at $p$.\par
As for $C_1$, we have $N_{C_1/\P^n}=(n-1)\O(n+2)$.  Then after blowing up the transverse
$\P_1^{e-n}=\P^{n-(2n-e)}$ to $C_1$ we get
\[N_{C_1'/X_1}\simeq (e-n)\O(n+2)\oplus (2n-e-1)\O(n+1).\]
This bundle is balanced and its upper subspace at $p$
comes from the first summand.
I claim that under the identification of the normal space with $T_pZ$, the upper
subspace has trivial intersection with $T_p\P_1^{e-n}$. To see this consider the projection
$X_1\to\P^{2n-e-1}$ which identifies $X_1=\P_{\P^{2n-e-1}}(\O(1)\oplus (e-n+1)\O)$ 
and maps $C_1$ to a general rational curve of degree $n-1$ in $\P^{2n-e-1}$.
Then (see \S \ref{vertical-normal}) 
the vertical (i.e. part killed by the blowdown map $X_1\to\P^{2n-e-1}$) 
subbundle of the normal bundle to $C_1$ has the form $K^*(1)$
where $K$ is the relative tautological subbundle which fits in an exact sequence
on $C_1$:
\[\exseq{K}{(e-n+1)\O\oplus\O(-n+1)}{\O(1)}.\]
Because the exceptional divisor $\P((e-n+1)\O)$ meets $C_1$ at $p$ only,
it follows that the induced map $(e-n+1)\O\to\O(1)$ vanishes at $p$
hence factors through an inclusion $\O\to\O(1)$. This yields a (locally split) inclusion
$(e-n)\O\to K$ and it follows easily that
	 $K\simeq (e-n)\O\oplus\oplus(-n+1)$ so
	\[K^*(1)\simeq (e-n)\O(1)\oplus \O(n).\] Now a nontrivial intersection of
	the upper subspace above with $T_p\P^{e-n}$ would yield a subbundle of $K^*(1)$
	which is a sum of copies of $\O(n+2)$ which is evidently impossible, proving our claim.

% kernel of $T_pZ\to T_p\P_1^{e-n}$.
 %Clearly 
 Hence, the two upper subspaces of $N_{C'_1/X_1}$ and $N_{C'_2/X_2}$ at $p$,
 both of which may be considered as subspaces of $N_{C_0/X_0}$ at $p$,  are transverse as such. 
Therefore $C_0$, which  is 
a locally complete intersection on $X_0$, hence also on $\cX$,
has normal bundle $N_{C_0/X_0}$ that  is balanced, positive and has $H^1=0$. 
Therefore as in \S \ref{subsub-sec}, $C_0$ smooths out to a smooth irreducible rational 
curve of degree $e$ in a nearby fibre
$\P^n\times t, t\neq 0$ with balanced normal bundle.\par
\ul{Case 3}: $e\geq 2n$.\par
This case is analogous to Case 2 except that we take $C_2\subset\P^n_2$ to be a general
(hence now nondegenerate)
curve of degree $e-n+1$ with $C_1\subset\P^n_1$ still a rational normal curve. 
By induction we may assume $N_{C_2/\P^n_2}$
	is balanced so we have
\[N_{C_2/\P^n_2}\simeq r\O(a^+)\oplus (n-1-r) \O(a^+-1),\]\[ 0<r\leq n-1, a^+=
\lceil((e-n+1)(n+1)-2)/(n-1)\rceil.\]
Then let $X_2$ be the blowup of $\P^n_2$ in a general  $\P^{n-r}_2$  meeting $C_2$
transversely in 1 point
 and $X_1$
be the  blowup of $\P^n_1$ in a general $\P^{r-1}_1$ meeting $C_1$ transversely
in 1 point, with $C'_1, C'_2$ being the birational transforms 
of $C_1, C_2$ resp.  So as above we have 
\[N_{C'_1/X_1}=(n-1-r)\O(n+2)\oplus r\O(n+1)\]
while 
\[N_{C'_2/X_2}=\O(a^+)\oplus (n-2)\O(a^+-1).\] Note that $a^+>n+2$ thanks to $e-n+1>n$, so
$N_{C'_1/X_1}$ cannot have an $\O(a^+)$ summand.
Hence we can argue as above that the upper subspaces of $N_{C'_1/X_1}$ and $N_{C'_2/X_2}$
at $p$ must have trivial intersection.

%****** 
%As to the upper subspace of $N_{C_2'/X_2}$ there are 2 cases.\par
%Case i:    $a^++2n-e-1\geq n-1$. Then the upper subspace will be a subspace of the 
%vertical subspace, kernel of $T_pZ\to T_p\P^{r}_2$, 
%in which case it will be automatically transverse to the upper
%subspace of $N_{C'_1/X_1}$.\par
%Case ii: $a^++2n-e-1\leq n-1$. Then the upper subspace it will be the
%sum of the same vertical subspace and the intersection 
%of the upper subspace of $(N_{C_2/{\P_2^n}})_p$
%with the (isomorphic) image of $T_p\P_2^{r}$ in $(N_{C_2/\P_2^n})_p$.
%Since the upper subspace of $N_{C_1'/X_1}$  coincides with $T_p\P_2^{r}$,
%the two upper subspaces are automatically transverse . This shows as above
%that $N_{C_1'\cup C'2/X_1\cup X_2}$
% is balanced so we can conclude as above.
\end{proof}

 The method of proof implicitly uses the notion 
of fang which will be revisited more explicitly in \S \ref{fangs}.

%In some ways, this proof is harder than that of the main theorems because there
%are proper upper subspaces on both sides while in the main theorems the upper
%subspace on one side is the whole space.
\section{Case $d=n$}\label{fano}
Our result is the following.
\begin{thm}\label{d=n-thm}
	Let $X$ be a general hypersurface of degree $n$ in $\P^n$, $n\geq 4$. Then 
	for any $e\geq n-1$, $X$ contains a nonsingular irreducible
	balanced rational curve of degree $e$.
	\end{thm}
\begin{cor}\label{ratcon-d=n-cor}
	Notations as above, $X$ is 
	separably $(\lfloor\frac{e-2}{n-2}\rfloor+1, e)$-
	rationally connected.
	\end{cor}
\begin{proof}[Proof of Corollary]
	Standard. Let $\cC/B$ be a component of the universal degree-$e$ rational curve in 
	$X$ containing a good curve as above and $\cC^q/B$ its $q$-th fibre power, which admits
	an obvious map
	\[f_q:\cC^q_B\to X^q\]
%	(note this $q$ was denoted $q+1$ in the introduction). 
	For $z=(C, p_1,...,p_q)\in
	\cC^q/B$, there is a derivative map
	\[df_q:T_z(\cC^q/B)\to \bigoplus T_{p_i}X\]
	taking the vertical part of the tangent space to $\bigoplus T_{p_i}C$,
	hence inducing $T_{[C]}B\to \bigoplus N_{p_i,C/X}$
	which is none other than the evaluation map
	\[H^0(N_{C/X})\to \bigoplus N_{p_i,C/X},\]
	with cokernel $H^1(N_{C/X}(-q))$.
	For $q=\lfloor\frac{e-2}{n-2}\rfloor+1=q_{\max}(e)$,
	the latter map is surjective by an evident $H^1$ vanishing, hence
	so is $f_q$ locally.
	\end{proof}
\begin{cor}
	Notations as above, for all $q\geq 2$, the minimal degree of a rational 
	curve in $X$ through $q$ general points is the expected 
	one, viz. $(q-1)(n-2)+2$, and the locus of such curves is reduced
	and  of the expected dimension.
	\end{cor}
\begin{proof}
	By Corollary \ref{ratcon-d=n-cor}, there exists a rational curve of degree
	$e=(q-1)(n-2)+2=e_{\max}(q)$ through $q$ general points. Because $q_{\max}(e-1)<e$,
	this $e$ is smallest.
	\end{proof}
\begin{proof}[Proof of Theorem \ref{d=n-thm}]
	
	Consider a 2-fan $P=P_1\cup_E P_2$ as in \S \ref{fans},
	so $P_1=B_p\P^n, P_2=\P^n, E=\P^{n-1}$. 
	%We will consider a degenerate form of $X$ in the the form of a fan
	On $P$, consider a general hypersurface $X_0$ of type $(n, n-1)$ as in \cite{hypersurf}, i.e.
	\[X_0=X_1\cup_F X_2,\]
	with $X_2$ is a general
	hypersurface of degree $n-1$ in $\P^n$ and $X_1=B_p\bar X_1$ is the blowup at $p$
	of a general quasi-cone $\bar X_1$ with quasi-vertex $p$, i.e. a hypersurface in $\P^n$ 
	of degree $n$ and multiplicity $n-1$
	at $p$, with exceptional divisor $F$ which
	is a degree-$n-1$ hypersurface in $E=\P^{n-1}$. Using projection from $p$,
	which maps $\bar X_1$ birationally onto $\P^{n-1}$,
	$X_1$ 
	can also be realized as
	the blowup of $\P^{n-1}$ in a general $(n-1,n)$
	complete intersection 
	\[Y=F_{n-1}\cap F_n, \deg(F_i)=i.\] 
	and in this realization $F$ is the birational transform of $F_{n-1}$
	and is isomorphic to the latter, while $Y$ can be identified as the set of lines through
	$p$ contained in $\bar X_1$.
	There is a family $\cX/\A^1\subset\mathcal P/\A^1$ 
	with general fibre $X\subset\P^n$ a general hypersurface of degree $n$
	and special fibre $X_0\subset P$. Concretely, if $p$ is the point $[1,0,...,0]$
	and we write the general degree-$n$ polynomial as $\sum\limits_{i=0}^n f_ix_0^{n-i}$
	with $f_i$ of degree $i$ in $x_1,...,x_n$
	then the equation for the quasi-cone $\bar X_1$ is $f_n+x_0f_{n-1}$, that of $X_2$
	is $\sum\limits_{i=0}^{n-1} f_ix_0^{n-i-1}$, that of $F$ and $F_{n-1}$
	is $f_{n-1}$, and that of $F_n$ is $f_n$.  \par
	\ul{Case 1}: $e\geq (n-1)^2$.\par
	Write $e=e_1n-a$ with \[e_1\geq n-1, a\leq n-1.\] 
	To construct a suitable curve in $X_0$ we proceed as follows.
	Let $C$ be a general rational curve of degree $e_1$
	in $\P^{n-1}$. Let $F_{n-1}\subset\P^{n-1}$ be a general
	hypersurface meeting $C$ transversely in $e_1(n-1)$ points.
	Note that $C\cap F_{n-1}$ is in (linearly) general position
	(i.e. any $n$ points linearly independent). Choose a subset
	$A$ from it with $|A|=a$ which we may assume consists
	of coordinate vertices $p_1,...,p_a$ (recall that $a\leq n-1$).
	We also assume $x_0,...,x_{n-1}$ are standard coordinates on $\P^{n-1}$
	(compatible with the $p_i$).\par
	Now I claim that that $\I=\I_{p_1,...,p_a/\P^{n-1}}(n-1)$, 
	in fact already $\I_{p_1,...,p_a/\P^{n-1}}(2)$, is globally generated.
	If $p$ is not in the span $S=\carets{p_1,...,p_a}$ then a general hyperplane containing 
	$S$ generates
	$\I$ at $p$. On the other hand if $p\in\carets{p_j:j\in J}$
	for some $J\subset \{1,...,a\}, |J|\geq 1$, with $J$ minimal, 
	then we can pick some $j\in J$ and a general hyperplane
	containing all $p_k, k\neq j$ plus another general hyperplane containing only $p_j$ to
	generate $\I$ at $p$.

%	The standard exact sequence
%	\[\exseq{\I_{C/\P^{n-1}}}{\I_{p_1,...,p_a/\P^{n-1}}}{\O_C(-p_1-...-p_a)}\] 
%	together with the fact that $\I_{C/\P^{n-1}}(n-1)$ and $\O_C(-p_1-...-p_a)(n-1)$
%	are globally generated and $H^1(\I_{C/\P^{n-1}}(n-1))=0$ 
%	shows that $\I_{p_1,...,p_a/\P^{n-1}}(n-1)$ is globally generated so $F_{n-1}$
	Now as $\I$ is globally generated, we may assume $F_{n-1}$ will miss any given finite set of points.
	I claim we can find a hypersurface $F_n$ through
	$A$ and no other points of $C\cap F_{n-1}$ and with given normal 
	hyperplanes to $C$ at $A$, i.e. given image of $T_pF_n$ in
	$N_{C/\P^{n-1}, p}$  for all $ p\in A$. Indeed, a degree-$n$ form through $p_i$
	 has no $x_i^n$ term and its tangent at  $p_i$
	corresponds to a term $x_i^{n-1}g_i$ with $g_i$ linear in $x_j, j\neq i$.
The $g_i$ may be chosen independently of one another and, identifying the normal space
to $C\subset\P^{n-1}$ at $p_i$ with $T_{p_i}F_{n-1}$,  $g_i$ specifies a
hyperplane in that normal space. Choosing the $g_i$ generally with this property
	and setting $F_n=\sum x_i^{n-1}g_i$ yields the desired $F_n$; choosing $F_{n-1}$
	general enough with given $p_1,...,p_a$ ensures that $F_{n-1}\cap F_n\cap C=\{p_1,...,p_a\}$.
\par	Now blow up \[Y=F_{n-1}\cap F_n\subset\P^{n-1}\] 
	to get $X_1$ and let $C_1\subset X_1$ be
	the birational transform of $C$. Because
	$C_1$ has balanced normal bundle and  $Y$ has general tangents at $Y\cap C$,
	$C_1\subset X_1$ also has balanced normal bundle, and it meets
	$F$ transversely in $e_1(n-1)-a$ points.
\par	 
Now let $C_2\subset X_2$ be $e_1(n-1)-a$ general lines so that
	$C_2\cap F=C_1\cap F$. As $X_2$ is a general hypersurface of degree $n-1$ it is easy to check that
	each of the lines has trivial 
	(i.e. globally free) normal bundle. Now in view of \S \ref{subsub-sec},
	 Corollary \ref{tree-cor}
	 applies
	and shows that $C_1\cup C_2$ smooths out to a smooth rational curve of degree $e$
	on a general hypersurface $X$ of degree $n$ in $\P^n$ with balanced normal bundle.
	
	\par
	\ul{Case 2}: $n-1\leq e<(n-1)^2$.\par
	We consider the same fan hypersurface $X_0=X_1\cup X_2$ as in Case 1 and in $X_0$ a curve 
	$C_1\cup C_2$ constructed as follows. Considering $X_1$ as $\P^{n-1}$ blown up in
	$Y=F_{n-1}\cap F_n$, we take for
	$C_1\subset X_1$	the birational transform of a general rational
	normal curve $C\subset\P^{n-1}$ meeting $Y$ in $a$ points with $a$ as above. 
	For $C_2$ we take $(n-1)^2-a$ general lines in $X_2$ corresponding to the
	points of $C\cap F_{n-1}$ not on $Y$. 
To conclude as in Case 1 that $C_1\cup C_2$ will smooth out to a balanced curve on $X$,
	it suffices to show that we can choose $C_1$ to have 
	balanced normal bundle in $X_1$. Let $2C$ be the first order neighborhood
	of $C$, with ideal sheaf $\I_C^2$, so we have an exact sequence
	\[\exseq{\check N}{\O_{2C}}{\O_C}\]
	where $\check N$ is the dual to 
	\[N_{C/\P^{n-1}}=(n-2)\O_C(n+1)\]
	where $\O_C(\l)$ is the line bundle of degree $\l$ on $C\simeq\P^1$.
	Because $\O_C(H)=\O_C(n-1)$,
	this sequence shows that \eqspl{ci}{H^1(\O_{2C}(kH))=0, k\geq 2} where $H$ is a hyperplane.
	%so that $\O_C(H)=\O_C(n-1)$.
	Let $Z=C\cap F_{n-1}, 2Z=2C\cap F_{n-1}$. 
	Then $\I_{Z/C}=\O_C(-H)$, $\I_{2Z/2C}=\O_{2C}(-H)$, hence we have exact sequences
	\[\exseq{\I_{C/\P^{n-1}}}{\I_{Z/\P^{n-1}}}{\I_{Z/C}}\]
	\[\exseq{\I^2_{C/\P^{n-1}}}{\I_{2Z/\P^{n-1}}}{\I_{2Z/2C}}.\]
 A theorem of Rathmann ca. 1991 (see \cite{rathmann-eh-conj}, Prop.4.2
	or \cite{vermeire}) shows that
	\eqspl{jurgen}{H^1(\I^2_{C/\P^{n-1}}(kH))=0, k\geq 3.} 
	Putting \eqref{ci} and \eqref{jurgen} we conclude
	 \[H^1(\I_{2Z/\P^{n-1}}(kH))=0, k\geq n-1\] since $n-1\geq 3$.
	This shows that the map $\rho: H^0(\O_{\P^{n-1}}(kH))\to H^0(\O_{2Z}(kH))$
	is surjective. Now the image of $T_pF_n$ in the normal space to $C$ at $p\in C\cap F_n$ 
	is just the hyperplane in the normal space, or equivalently, element in the conormal space,
	which corresponds to $\rho(F_n)$. By surjectivity of $\rho$, therefore, $F_n$ can be chosen 
	to have general tangent hyperplanes 
	at the points $p\in Y\cap C$ modulo $T_pC$, which makes $N_{C_1/X_1}$ balanced. 
	Therefore again $C_1\cup C_2$
	smooths out to a balanced rational curve on $X$.
		\end{proof}
%\begin{rem}
%	The proof does not extend to the case of hypersurfaces of degree $d<n$
%	because in that case lines on $F_{d-1}$ do not have 
%	trivial normal bundle and, especially if there are many lines involved, it
%	is not clear that the gluing of normal bundles of $C_1$ and $C_2$
%	can be chosen general enough so as to be able to apply Corollary \ref{tree-cor}.
%	However the analogous argument does yield the following statement:
%	for $d<n$ and a general hypersurface $X$ of degree $d$ in $\P^n$, 
%	the blowup of $X$ in a general linear space section of codimension
%	$n+1-d$ contains a balanced rational curve of every degree $e\geq n-1$.
%	This works because in the appropriate blowup of $F_{d-1}$, lines have
%	trivial normal bundle.
%	\end{rem}

%\section{Fangs and the case $d<n$}\label{fang}
\section{More preliminaries}\label{prelim2}
	In the next section we will construct, for infinitely many degrees $e$,
	rational curves of degree $e$ with balanced normal bundle on
	a general hypersurface $X$ of degree
	$d$ in $\P^n$ with $3\leq d\leq n-1$. 
	Here we collect a few more
	preliminary results and constructions that we will need.

	%\section*{Preliminaries}
	\subsection{Curves in projective bundles}\label{vertical-normal}
	Let $G$ be a vector bundle on a variety $B$, with associated projective
	(quotient) bundle $\P(G)\stackrel{\pi}{\to} B$. 
	Given a parametrized curve $c_0: C_0\to B$, a lifting  of 
	$c_0$ to $c:C_0\to \P(G)$ corresponds to an invertible quotient
	$G_0:=c_0^*G\to L_0$. In this case we have
	\[c^*\O_{\P(G)}(1)=L_0.\]
	There is an exact sequence for the normal bundle
	\[\exseq{N_{C_0/\P(G)/B}}{N_{C_0/\P(G)}}{c^*N_{C_0/B}}\]
			and, setting $K_0:=\ker(G_0\to L_0)$ (the relative tautological subbundle),
	the \emph{vertical} normal bundle, i.e. $N_{C_0/\P(G)/B}$, is given by
		\[N_{C_0/\P(G_0)/B}=K^*_0\otimes L_0 .\]
		The \emph{horizontal} normal bundle is by definition $c^*N_{C_0/B}$.
	\subsection{Some blowups}\label{blowup}
		%\vskip 1cm
	%
	%/*********************
	
	%
	%Let $\Gamma_1=\P^{n-m-1}, \Gamma_2=\P^m$ be disjoint linear spaces
	%in $\P^n$. Let $b_i:Z_i\to\P^n$ be the blowup of $\Gamma_i, i=1,2$, and
	%$b:Z\to\P^n$ the blowup of $\Gamma_1\coprod\Gamma_2$. Thus we have 
	%a diagram
	%\eqspl{}{
	%\begin{matrix}
	%	&&Z&&\\
	%	&\swarrow&&\searrow&\\
	%	Z_1&&&\  \  \  \  \  \  \  Z_2\\
	%	&\searrow&&\swarrow&&\\
	%	&&\P^n&&
	%	\end{matrix}
	%	}
	%Now there are projecrtion maps
	%\[\pi_1:Z_1\to\P^m, \pi_2:Z_2\to\P^{n-m-1}, 
	%\pi:Z\to\P^m\times\P^{n-m-1}\]
	%that fit in a diagram
	%\[
	%%\eqspl{}{
	%\begin{matrix}
	%	Z_1&\leftarrow&  Z&\to& Z_2\\
	%	\downarrow&&  \downarrow&&\downarrow\\
	%	\P^m&\leftarrow&\  \  \  \P^m\times\P^{n-m-1}&\to&\P^{n-m-1}
	%	\end{matrix}	
	%%}
	%\]
	%In fact, $Z_1=\P_{\P^m}(1, 0^{n-m}):=\P(\O(1)\oplus(n-m)\O)$,
	%$Z_2=\P_{\P^{n-m-1}}(1, 0^{m+1})$, while
	%$Z=\P_{\P^m\times\P^{n-m-1}}(p_1^*\O(1)\oplus p_2^*\O(1))$. Thus $Z$ 
	%contains two disjoint sctions $S_1, S_2$, which project respectively to
	%the exceptional divisors $E_1=\P(0^{n-m})=\P^m\times\P^{n-m-1}\subset Z_1$,
	%$E_2=\P(0^{m+1})=\P^{n-m-1}\times\P^m\subset Z_2$.\par
	%Now let $Y_0\subset\Gamma_1\times\Gamma_2$ be a hypersurface
	%of bidegree $(1, d-1)$ and let $X=\pi\inv(Y_0)\subset Z$, and
	%let $X_1, X_2$ be the images of $X$ in $Z_1, Z_2$. Note that
	%$X_1$ is a projective subbundle of $Z_1$, namely $\P(G)$
	%where $G$ is the cokernel of a suitable map $\O(-(d-1))\to\O(1, 0^{n-m})$
	%on $\P^m$. Likewise $X_2$ is fibred over $\Gamma_2$ with fibres
	%hypersurfaces of degree $d-1$ in $\P^{m+1}$. 
	%*******************/
This construction will be used to construct components of fangs and fang hypersurfaces.
	Let $b:Z\to\P^n$ denote the blow-up of $\P^n$ in $\P^{n-m-1}$ and
	$\pi:Z\to\P^m$ the morphism induced by linear projection with center
	$\P^{n-m-1}$, whose fibres are fibres $\P^{n-m}$.
	Note that via $\pi$, $Z$ is a projective bundle
	 \[Z=\P_{\P^m}(\O(1)\oplus(n-m)\O):=\P(1, 0^{n-m})\]
	and  that the exceptional divisor of $b$ is 
	\[E=\P_{\P^m}((n-m)\O)=\P^m\times\P^{n-m-1}.\]
	A hypersurface $X$ of type $(d,b)$ on $Z$ is an element of
	the linear system $|b^*\O(d)-bE|$ and $X$
	maps to a hypersurface of degree $d$ in $\P^n$ with multiplicity
	$b$ on $\P^{n-m-1}$ and $X$
	meets $E=\P^m\times\P^{m-n-1}$ in a divisor $Y$   of bidegree
	$(d-b, b)$ on $E$. 
	The fibres of $\pi|_X$ are hypersurfaces of degree $d-b$
	in $\P^{n-m-1}$.\par
	Here are the two case of the construction to be used below.
\subsubsection{Special Case 1: $b=d-1$}\label{case1}  This means $X$ maps to a quasi-cone in $\P^n$. 
	 Then $X$ is a projective subbundle
	of $Z$ of the form $\P(G)$ where $G$ fits in an exact sequence
	\eqspl{G}{\exseq{\O_{\P^m}(-(d-1))}{\O_{\P^m}(1)\oplus(n-m)\O_{\P^m}}{G}.
	}
	In terms of coordinates, the image of $X$ in
	$\P^n$, for any $b$, has an equation of the form
	\eqspl{equation}{
		f=\sum\limits_{i=0}^{d-b}\sum\limits_ja_{d-i,j}(x)b_{i,j}(y)
	}
	where 
	$x_0,...,x_m, y_1,...,y_{n-m}$ are linear coordinates
	on $\P^n$ such that $x_0,...,x_m$ cut out $\P^{n-m-1}$,
	and each $a_{k,j}, b_{i,j}$ is homegenous of degree $k, i$ respectively.
	The term $i=d-b$ in the same equation, viewed as bihomogeneous form of bidegree $(b, d-b)$,
	also yields the equation of $Y$ in $E$.
	%and, when $e=d-1$, the left map in \eqref{G}.\par
	When $b=d-1$, we can write the equation $f$ in the form
	\eqspl{f}{f=a_d(x)+\sum\limits_{j=1}^{n-d+1}y_ja_{d-1, j}(x).}
	The left map in \eqref{f} is just given by $(a_d, a_{d-1, 1},...,a_{d-1, n-d+1})$.
\subsubsection{Special Case 2: $b=1, d=n-m$.}\label{case2}	
	Then the general fibre of $\pi_X$ is a  hypersurface
	of degree $n-m-1$ in $\P^{n-m}$. As is well known, a
	general such hypersurface is filled up by lines $L$ with trivial
	normal bundles (cf. Remark \ref{trivial} below), 
	and of course $L$ must meet $\P^{n-m-1}$.
	Now the birational transform of $L$ in $X$ has normal bundle which
	is an extension of trivial bundles, hence trivial as well.
	Consequently:
	\begin{lem}\label{x2}
		Notations as above, let $X$ be a general hypersurface of type $(n-m, 1)$ in $Z$.
%		and assuming general choices and $m=n-d$,
Then		there is a filling family in $X$ of birational transforms of lines
		meeting the projection center $\P^{n-m-1}$, hence contained in fibres
		of $\pi$,
	whose general member has trivial normal bundle in $X$.	
	\end{lem}
	\subsection{Fans generalized}\label{fangs}
	This is a slight generalization of the fans considered  above
	and in \cite{hypersurf}, and has already occurred in \S \ref{pn}.
	Let $\pi_1:Z_1=B_{\P^{n-m-1}}\P^n\to\P^m$ be as in \S\ref{blowup}
	with blowdown map $b_1:Z_1\to\P^n$ and exceptional divisor $E_1$. 
	Let $\pi_2:Z_2=B_{\P^m}\P^n\to\P^{n-m-1}$ be the analogous
	object, based on blowing up $\P^m$, with blowdown map $b_2$
	and exceptional divisor $E_2$. Note that both $Z_1$ and
	$Z_2$ have exceptional divisor $E=\P^m\times\P^{n-m-1}\simeq E_1\simeq E_2$. 
	The normal-crossing
	variety
	\[Z_0=Z_1\cup_EZ_2\]
	is called a \emph{generalized fan} or \emph{fang} of type $(n,m)$.\par
	A flat morphism $\cZ\to B$ is called a \emph{relative fang} of
	type $(n,m)$ if each fibre is either $\P^n$ or a fang of type $(n,m)$.
	A standard way to construct a relative fang is to blow up
	the subvariety $\P^{n-m-1}\times 0$ in $\P^n\times\A^1$.
	There $Z_1$ and $Z_2$ are, respectively, the birational transform of $\P^n\times 0$
	and the exceptional divisor, and both have projective bundle structure:
	\[Z_1=\P(G_1), Z_2=\P(G_2)\]
	where
	\[G_1=\O_{\P^m}(1)\oplus (n-m)\O_{\P^m}, G_2=\O(1)_{\P^{n-m-1}}\oplus (m+1)\O_{\P^{n-m-1}}.\]
	We will denote the latter as $Z_1=\P_{\P^m}(1,0^{n-m}), Z_2=\P_{\P^n-m-1}(1, 0^{m+1})$.
Also each $Z_i$ is endowed with the $\O(1)$ induced by from $G_i$ which is
	also the pullback of $\O_{\P^n}(1)$ by the blowdown map $b_i:Z_i\to\P^n$.
	We denote the exceptional divisor of $b_i$ by $E_i$, and we also have
	\[E_1=\P((n-m)\O_{\P^m}), E_2=\P((m+1)\O_{\P^{n-m-1}}),\]
	hence of course
	\[E=E_1\simeq E_2\simeq \P^m\times\P^{n-m-1}.\]\par
The case $m=0$ or $m-n-1$ reduces to the fan construction used previously,
	Now assume $0<m<n-1$ and $e<d$. Then
	the linear system $|dH-eZ_2|$ on $\cZ$, 
	where $H$ is the pullback of a hyperplane in
	$\P^n$, restricts as follows.
	\begin{itemize}\item on the general fibre, to $|dH|$;\item
	on $Z_1$,  to $|db_1^*H_{\P^n}-eE_1|$, i.e. the birational transform on $Z_1$
	of the system of hypersurfaces of degree $d$ on $\P^n$ 
	with multiplicity $e$ on $\P^{n-m-1}$;\item
	on $Z_2$ to $|db_2^*H_{\P^n}-(d-e)E_2|$;\item on 
	$E=\P^m\times\P^{n-m-1}$ to the linear system of hypersurfaces
	bidegree $(e, d-e)$.
	\end{itemize}
\begin{example}\label{curve-degree} Taking $d=1, e=0$  shows that a curve of degree $k$ in 
$\P^n$ can specialize to a curve $C_1\cup C_2\subset Z_0$ where $C_1\subset Z_1$ maps to
a curve $b_1(C_1)$ of degree $k_1$ in $\P^n$ while $C_2$ 
a curve
$b_2(C_2)\subset\P^n$ of degree $k-k_1+\l$ meeting te blown-up $\P^m$ in $\l$ points (hence $C_2$
maps to a curve $\pi_2(C_2)$ of degree
$k-k_1$ in the base $\P^{n-m-1}$). 
 Here $\l=C_1.E=C_2.E$ and $C_1\cap E=C_2\cap E$ consists of $\l$ points.
 \par Similarly for the dual case $d=e=1$.
	\end{example}
	The foregoing construction may obviously be extended to the case of
	more than 2 components but we don't need this.
	
	\subsection{Balanced extensions and kernels}
	An extension of balanced vector bundles is balanced when their slopes
	are roughly equal. This is useful for constructing balanced bundles.
	\begin{lem}\label{extension}
		Let
		\[\exseq{E_1}{E}{E_2}\]
		be an exact sequence of  vector bundles on $\P^1$,
		of respective slopes $s_1, s, s_2$. Assume $E_1, E_2$ are 
		balanced and 
		\eqspl{extension-eq}{
		\lfloor{s_1}\rfloor=\lfloor{s_2}\rfloor.}
	
		Then $E$ is balanced and $\lfloor{s}\rfloor=\lfloor{s_1}\rfloor$.
		Moreover the extension splits.
	\end{lem}
	The proof may be left to the reader.\par
	Balancedness is also inherited by the kernel of a general map to a vector bundle:
		\begin{lem}\label{kernel}
		Let $E$ be a balanced bundle on $\P^1$ and $\phi: E\to L$ a 
		sufficiently general surjection
		to a vector bundle. Then $\ker(\phi)$ is balanced.
	\end{lem}
	\begin{proof} By an obvious induction we may  assume $L$ has rank 1.
	%	hence that $L=\O$. Write $E=(r-r^-)\O(a^-+1)\oplus r^-\O(a^-), r^->0$.
	%	If $a^->0$, the only map $E\to L$ is zero, whose kernel is balanced.
	%	If 
		Because balancedness is open it suffices to prove: given $E$ 
		balanced of slope $s$ and an integer $\l\geq \lfloor s\rfloor$, there
		exists a balanced bundle $K$
		with $c_1(K)=c_1(E)-\l, \rk(K)=\rk(E)-1$ and a locally split injection $K\to E$.
		We may assume $E=r_+\O(1)\oplus (r-r_+)\O$, so $\l\geq 0$.
		Write
		\[\l=q(r-1)+p, 0\leq p<r-1.\]
		If $p\leq r_+$ we can take
		\[K=((r_+-p)\O(1-q)\oplus (r-r_++p)\O(-q).\]
		If $p>r_+$ we can take
		\[K=(r-p+r_+)\O(-q)\oplus (p-r_+)\O(-q-1).\]
		Clearly, a general map $K\to E$ is locally split injective.
	\end{proof}
	\begin{rem}\label{trivial}
		The Lemma yields a quick proof of the fact that a general line $L$ on 
		a general hypersurface $X$ has balanced normal bundle: 
		indeed if $x_2,...,x_n$ are linear equations for $L$ and $\sum\limits_{i=2}^nx_if_i$
		is an equation for $X$ then $N_{L/X}$ is the kernel of the general map
		$N_{L/P^n}\stackrel{(x_2,...,x_n)}{\to}\O_L(d)$.\par
		In particular,
		a general line on a general hypersurface of degree $n-1$ in $\P^n$
		has trivial normal bundle.
	\end{rem}
The following is a close analogue of Lemma \ref{mod}, showing that
a suitable elementary down modification brings a bundle closer to balance.
\begin{lem}\label{straightening}
	Let 
	\eqspl{e}
	{E=r_0\O(a)\oplus r_1\O(a-1)\oplus...\oplus r_b\O(a-b),\ \ r_0>0, r_1,...,r_b\geq0}
be a bundle on $\P^1$ and for $p\in\P^1$ denote by $\kk_p$ the  skyscraper sheaf $\kk(p)$ at $p$.
Then there exists a map $\phi:E\to\kk_p$
such that $E_1=\ker(\phi)$ has the form
\[E_1=(r_0-1)\O(a)\oplus (r_1+1)\O(a-1)\oplus...\oplus  r_b\O(a-b).\]
	\end{lem}
The proof is obvious as it suffices to take a map that is nontrivial
on one $\O(a)$ summand and zero on all other summands. Then applying the Lemma
$r_0$ times we conclude that there is a map $E\to\bigoplus\limits_{i=1}^{r_0}\kk_{p_i}$
with kernel
\[E_{r_0}=(r_0+r_1)\O(a-1)\oplus...\oplus r_b\O(a-b).\]
Thus passing from $E$ to $E_{r_0}$ decreases by 1 the degree difference between the
most positive subbundle and least positive quotient bundle.
Continuing in this manner at least $\rk(E)-r_b$ many times, 
and using openness of balancedness, we conclude
\begin{cor}\label{straightening-cor}
	(i) Notations as above, the kernel of a sufficiently general map 
	$E\to\bigoplus\limits_{i=1}^{s}\kk_{p_i}$ is balanced (resp. perfectly balanced) 
	provided $s\geq \sum\limits_{i=0}^{b-1}(b-i)r_i$
	(resp. $s=\sum\limits_{i=0}^{b-1}(b-i)r_i$).
	\par (ii) In particular, if $E$ is balanced with upper rank $r^+$
	then a sufficiently general map $E\to\bigoplus\limits_{i=1}^{r^+}\kk_{p_i}$
	has perfectly balanced kernel.
	\end{cor}
\begin{rem}[non-essential] 
	We may define the 'unbalanced degree' of $E$ as above as $u(E)=\sum r_i(b-i)$
	and its 'lower rank' $r^-(E)$ as the rank of its lowest-slope quotient, i.e. $r_b$.
	Clearly $u(E)\geq \rk(E)-r^-(E)$ with equality iff $E$ is balanced. 
	Then Lemma \ref{straightening} implies that, unless $E$ is perfectly balanced,
	we have $u(E_1)=u(E)-1$. Also, $r^-(E_1)=r^-(E)$ unless $r^-(E)=\rk(E)-1$.
	Since $u(E)\geq 0$, the modification $E\mapsto E_1$ yields a perfectly
	balanced bundle after $u(E)$ many iterations.
	This yields another way
	to deduce Corollary \ref{straightening-cor}.
	\end{rem}
	\section{Fangs and the case $d<n$}\label{d<n}
	Here we will prove the results on balanced rational curves
	on a general hypersurface $X$ degree $d<n$.
	Since it is known by 
	the result of Riedl-Yang that the family
	of rational curves of degree $e$ on $X$ is irreducible at least if $d<n-1$, it follows that almost all
	of these curves have balanced normal bundle. This 
	construction has 
	consequences for minimal rational connectivity. \par
	To state the results we need a definition. Fix n, d with $3\leq d\leq n-1$.
	An integer $e$ is said to be \emph{accessible via} $e_0$
	 if   $e\geq e_0\geq d-1$ and
	\eqspl{accessible}{
		\lfloor \frac{-de_0+e}{n-d}\rfloor+e=e_0+\lfloor\frac{2e_0-2}{d-2}\rfloor.	
	} $e$ is accessible if it is accessible via some $e_0$.
Also recall (see \S \ref{numerology-sec}) that an integer $e$ is said to be point-minimal
if $q_{\max}(e-1)<q_{\max}(e)$.
	\begin{thm}\label{d<n-thm}
	(i) A general hypersurface
	$X$ of  degree $d<n$
	in $\P^n$ contains balanced rational curves of 
	any accessible degree $e$ and is 
	separably rationally 
	$(q_{\max}(e)-1)$- connected for any accessible $e$.
	\par (ii) For  $2<d<n$, the the set of accessible $e$ contains the intersection
	of $[d-1, \infty)$ with $a(d, n-d)$ many
	%	$(n-d)(d-5/2)$ 
	congruence classes mod $d(n-2)$ where $a(n,n-d)=(n-d)d-\mathrm{(linear\ terms)}$. \par
	(iii) For $3<d<n-1$, the
	set of $e$ which are both accessible and point-minimal
	contains the intersection of $[d-1, \infty)$ with $(n+1-d)/2$ many congruence 
	classes mod $d(n-2)$ and the set of interpolating $q$
	%and $X$ is minimally separably rationally $q$-connected  
	%for all $q$
	contains some ray intersected with $(n+1-d)/2$ many congruence classes mod $(n+1-d)d$.
\end{thm}
Below we prove part (i). Parts (ii), (iii) follow from this by arithmetic:
 see Example \ref{d=n-1-example}
(for $d=n-1$) or the Appendix by M. C. Chang.
	\begin{example}\label{d=n-1-example}
	Take $d=n-1, n\geq 4$. Write $e_0=k(n-3)+r, k\geq 1, 0\leq r<n-3$.
	Then either 
	\[n\ \ \mathrm{even,}\ \  n\geq 6, \ 0<r\leq \frac{n-3}{2},\   
	e=\binom{n-1}{2}k+\frac{nr}{2}\] 
	or 
	\[n, r \mathrm{\ \ both\  odd,\ \ } n\geq 5,\  r\geq\frac{n-1}{2},
	e=\binom{n-1}{2}k+\frac{nr+1}{2},\]
	or
	\[ n\mathrm{\ \ odd},\  r\mathrm{\ \ even},\  r\leq \frac{n-3}{2},
	\ e=\binom{n-1}{2}k+\frac{r}{2}, \]
	or else
	\[n=4, r=0, e=3k-1.\]
	
	Thus, the accessible degrees cover about  $(n-3)/2$ 
	of the possible congruence classes of $e$ mod $\binom{n-1}{2}$.\par
	The point-minimal condition on the accessible $e$ is that
	the remainder of $2e-2$ mod $n-2$ should equal 0 or 1. 
	Considering those $e$ that are both accessible and point-minimal
	yields that $X$ is minimally
	rationally $(q-1)$-connected for $q-1=(n-1)k+1$ if  $n$ is even 
	(resp. $q-1=(n-1)k+2$ if  $n$ is odd),
	for any $k\geq 1$.\qed
\end{example}

\begin{proof}[Proof of Theorem \ref{d<n-thm}, (i)] The proof is based on a relative fang
	as in \S \ref{fangs}. Thus,
	fixing integers $d<n$, let $\cZ\to\A^1$ be a relative fang
	of type $(n, m), m=d-1$, with general fibre $\P^n$ and
	 special fibre $Z_0=Z_1\cup Z_2$. 
	Thus
	\[Z_1=\P_{\P^m}(1, 0^{n-m}), Z_2=\P_{\P^{n-m-1}}(1, 0^{m+1}).\]
	Consider a general member of the linear system $|dH-(d-1)Z_2|$
	on $\cZ$, whose general fibre over $\A^1$ is a general hypersurface
	of degree $d$ in $\P^n$, and let
	\[X_0=X_1\cup X_2\]
	be its special fibre over $\A^1$. Thus, $X_1=\P(G)\to\P^m$ as in \S \ref{case1}, while
	$X_2$ fibres over $\P^{n-m-1}$ with general fibre a general hypersurface of degree
	$d-1=(m+1)-1$ in $\P^{m+1}$ \S \ref{case2} (beware the switch in notation,
	interchanging $n-m$ and $m+1$). The idea, as in the proof of Theorem \ref{d=n-thm},
	is to construct a good curve $C_1\cup C_2\subset X_0$
	which will smooth out to a balanced curve on $X$.
	In fact, $C_1\subset X_1$ will be balanced and $C_2\subset X_2$ will be perfectly balanced.\par
We first construct $C_1$. To this end	consider a general rational curve  $C_0$ of
	degree $e_0$ in $\P^m$, and let $C_1\subset X_1$ be a general 
	degree-$e$ lifting of $C_0$,
	which corresponds to a general surjection
	\eqspl{surjection}{\psi: G_{C_0}\to\O_{C_0}(e)}
	Let $K=\ker(\psi)$, i.e. the restriction of the relative tautological
	subbundle of the projective bundle $\P(G)$.
	Then $C_1$ meets the exceptional divisor $E=\P^m\times\P^{n-m-1}, m=d-1$,
	in $e-e_0$ points. By Proposition \ref{pn-prop}, as soon as $e_0\geq m$, 
	the normal bundle $N_0=N_{C_0/\P^m}$
	is balanced, of slope $ s_0=\frac{(m+1)e_0-2}{m-1}$. 
	By \S \ref{vertical-normal}, the 'relative'
	or vertical normal bundle $N_{C/X_1/\P^m}$ is just
	$K^*(e)$. Thus,
	we have an exact sequence
	\eqspl{normal}{
		\exseq{K^*(e)}{N_{C_1/X_1}}{N_0}.
	}
	\begin{lem}\label{k-balanced}
		Notations as above, $K$ is balanced.
	\end{lem}
\begin{proof}
Assume first that $e_0\leq n-d+1$. Recall the surjection
\[s:\O(e_0)\oplus(n-m)\O\to G_{C_0}.\]	 
%************We have an exact diagram
%\eqspl{}{
%\begin{matrix}
%	&&&\O(e_0)&=&\O(e_0)&\\
%	&&&\downarrow&&\downarrow&\\
%0\to&\O(-(d-1)e_0)&\to&\O(e_0)\oplus(n-m)\O&\to&G_{C_0}&\to 0\\
%	\end{matrix}	
%}
%****/\par
By generality of $\psi$ we may assume the induced map
	$\O(e_0)\to\O(e)$, i.e. the restriction of $\psi\circ s$ on the $\O(e_0)$ summand, 
	is injective and general and let $\tau$ be its cokernel, which is a
	skyscraper of the form
$\tau=\bigoplus\limits_{i=1}^{e-e_0}\kk_{p_i}$,
where $p_1,...,p_{e-e_0}\in C_0$ are distinct points. Let $G_1$ be the cokernel of the map
	$\O(e_0)\to G_{C_0}$, so we have an exact sequence
			
		\eqspl{g1}{\exseq{\O(e_0)}{G_{C_0}}{G_1}.}
			
Then from \eqref{g1} we get another exact sequence on $C_0$:
%\[\exseq{\O(e_0)}{G_{C_0}}{(n-d+1)\O},\]
\eqspl{kg1}{\exseq{K}{G_1}{\tau}}
 %$\tau$ is the cokernel of $\O(e_0)\to\O(e)$, hence a 
%skyscraper sheaf of the form
 %$\bigoplus\limits_{i=1}^{e-e_0}\kk_{p_i}$,
where the map $\psi_1:G_1\to \tau$ is general. %Therefore $K$ is balanced.	

%/*************	
%This follows from
%Lemma \ref{kernel} provided we can prove that $G_{C_0}$ is balanced.
%To prove the latter assume first that $e_0\leq n-d+1$. Now as
%in \eqref{G} and \eqref{equation}, $G_{C_0}$ fits
Now I claim $G_1$ is balanced. To this end note $G_1$ fits in an exact sequence
\eqspl{GC}{\exseq{\O(-(d-1)e_0)}{ (n-d+1)\O}{G_{1}}}
where the left map is given by restriction of part of the equation
in $\P^n$ of the quasi-cone \eqref{f}, image of $X_1$, which is general as such.
Now by the maximal rank property of general rational curves \cite{ballico-ellia-max-rank},
the restriction maps $H^0(\O_{\P^m}(\l))\to H^0(\O_{C_0}(\l e_0))$
are surjective, $\l=d, d-1$, because
$de_0+1\leq d(n-d+1)+1\leq\binom{d+n}{n}$. Therefore the left map in
\eqref{GC} is general, so by Lemma \ref{kernel}
$G_1$ is balanced. Finally, since the map $\psi_1$ is general, it follows 
by Corollary \ref{straightening-cor}%Lemma \ref{torsion} 
that $K$ is balanced,
 provided $e_0\leq n-d+1$.\par
Now assume $e_0>n-d+1$ and write $e_0=e_1+(n-d+1)$ and consider a general
connected curve of the form $C_{11}\cup C_{12}$ where $C_{11}$ (resp. $C_{12}$) is
a general rational curve of degree
$e_1$ (resp. $n-d+1$). By induction, $G_{C_{11}}$ is balanced while $G_{C_{12}}$
is perfectly balanced (twist of a trivial bundle). Then using 
Theorem \ref{main-thm}, case (ii)', it follows that $C_{11}\cup C_{12}$ smooths to a curve $C_0\subset\P^m$ 
on which $G$ is balanced.
This proves Lemma \ref{k-balanced}.
%******/
\end{proof}

%to $G^*$ in the role of $K^*$
%	(cf. \eqref{G}))
%	shows that $G_{C_0}$ is balanced. Then applying the lemma again to 
%	the surjection $G_{C_0}\to \O(e)$ shows that $K$ is balanced.
	Now that we know $N_0$ and $K$ are balanced, we can apply
	 Lemma \ref{extension} to the exact sequence \eqref{normal},
	to conclude that $N_{C_1/X_1}$ is balanced provided the 
	matching condition \eqref{accessible} holds, i.e. provided $e$ is
	accessible via $e_0$. \par
	Note that the curves $C_1$ produced above meet $E$ 
		in $e-e_0$ points, say $p_1,...,p_{e-e_0}$. Let $L_i\subset X_2$
		be a line in $X_2$ through $p_i$ contained in the fibre through $p_i$
		of $X_2/\P^{n-m+1}$. Then by Lemma \ref{x2},  each $L_i$ has trivial normal bundle
		in the fibre, hence in $X_2$. Now let $C_2=\bigcup L_i$.
		%Now we argue as in \S \ref{fano} and attach lines in $X_2$ with trivial normal bundle at
		%those points, (cf. Lemma \ref{x2}), 
		Then $C_1\cup C_2$ is an lci curve in $X_0$ with 'balanced'
		normal bundle as in \S \ref{fano}, which smooths out to a balanced
		rational curve of degree $e$ on a general hypersurface of degree $d$. This proves
				Theorem \ref{d<n-thm}, Part (i).

%	This proves  Theorem \ref{d<n-thm}.
	\end{proof}
\begin{rem}
	The use of  the maximal rank property of the curve $C_0$ could be avoided 
	at the cost of forcing $e$ to satify a lower bound. 
	Because by \eqref{GC}, $G_1$ has no negative quotient
	hence its most positive line subbundle is $\O(a), a\leq e_0(d-1)$ and we can apply
	Lemma \ref{straightening} to conclude that $K$ will be balanced provided $e\geq de_0$.
	\end{rem}

		\begin{rem}
			Extending the fang method beyond the accessible values of $e$ would require
			taking $m\neq d-1$, and hence attaching lines with
			balanced, but nontrivial
			normal bundle. One would have to prove a general position property
			for the upper subspaces of these normal bundles, as in Theorem 
			\ref{main-thm}, (i). 
			This seems difficult.
			\end{rem}
		
%		\begin{rem}\label{mcc}
%			In the appendix below, M. C. Chang 
%			proves that  the set of accessible
%			values of $e$  is the intersection of $[d-1, \infty)$
%			with the union of at least $(n-d)(d-\frac{5}{2})$ many congruence classes 
%			$\mod (n-2)d$. 
%			 In particular, in the large-index case $n>> d>>0$, 
%		 the density of the set of $e$
%			is close to 1, while at the other extreme $d=n-2>>0$, the density is at least
%			(and actually approximately equal to)
%			$2/n$. The Appendix also gives bounds on the density of the 
%			set of point-minimal $e$ and corresponding $q$. The bounds are
%			probably not optimal.
%			\end{rem}

		%	which yields $e=k\binom{n-2}{2}+n/2$ so for
		%	each such $e$ we get a balanced rational curve of degree $e$
		%	on a general hypersurface of degree $n-1$ in $\P^n$.

\vfill\eject

\bibliographystyle{amsplain}
\bibliography{../mybib}
%\bibliography{../mybib}
\vfill\eject
%\documentclass[12pt]{article}
%    \usepackage{amsmath, amssymb, amsthm, amsxtra   }
%    %\numberwithin{equation}{section}
%\def\diam{\text{diam }}
%\def\supp{\text{supp }}
%\def\be{\begin{equation}}
%\def\ee{\end{equation}}
%\def\ve{\varepsilon}
%\def\vp{\varphi}
%\def\arrowk{^\to{\kern -6pt\topsmash k}}
%\def\arrowK{^{^\to}{\kern -9pt\topsmash K}}
%\def\arrowt{^\to{\kern -6pt\topsmash t}}
%\def\arrowr{^\to{\kern-6pt\topsmash r}}
%\def\arrowvp{^\to{\kern -8pt\topsmash\vp}}
%\def\tk{\tilde{\kern 1 pt\topsmash k}}
%\def\barm{\bar{\kern-.2pt\bar m}}
%\def\barN{\bar{\kern-1pt\bar N}}
%\def\barA{\, \bar{\kern-3pt \bar A}}
%\def\curl{\text { curl\,}}
%\def\G{\underset\cal G \to +}
%\def\div{\text{ div\,}}
%%\hsize = 6.2true in \vsize=8.2 true in
%\def\be{\begin{equation}}
%\def\ee{\end{equation}}
%\numberwithin{equation}{section}
%\begin{document}
    \theoremstyle{plain}
    \newtheorem{theorem}{Theorem}
    \newtheorem{lemma}{Lemma}
    \newtheorem{proposition}{Proposition}
    \newtheorem{corollary}[theorem]{Corollary}
    %only theorems and corollaries are on the same counter in this scheme
    \theoremstyle{definition}
    \newtheorem*{definition}{Definition}
    \theoremstyle{remark}
    \newtheorem*{remark}{Remark}

%\title{ \small{ }
%\footnote{2010 {\it Mathematics Subject
%Classification}.Primary 11B25.} \footnote{{\it Key words}. arithmetic progressions, quantitative %Nullstellensatz.}}
%\author{\small{Mei-Chu Chang}\footnote{Research partially
%financed by the NSF Grants~DMS~1600154.}\\ \texttt{\small{Department of Mathematics}}\\
%\texttt{\small{University of California, Riverside}}\\\texttt{\small
%mcc@math.ucr.edu}}

%$\qquad\qquad\qquad\qquad\qquad\qquad${\bf Appendix}
\section*{Appendix}
\medskip

$\qquad\qquad\qquad\qquad\qquad\qquad\qquad$ By Mei-Chu Chang\footnote
{Department of Mathematics,
	University of California, Riverside CA 92521,
 mcc@math.ucr.edu. 
 
 Research partially
financed by the NSF Grants~DMS~1764081.}

\bigskip

\bigskip

In this appendix, we prove Theorem A.1 and Theorem A.2 below.

\medskip

\noindent{\bf Theorem A.1.} {\it Let $2<d<n$ be integers. Then the set of integers $e$ which are accessible, i.e.  such that for some
$e_0$, $d-1\leq e_0\leq e$, one has
$$\bigg\lfloor \frac{-de_0+e}{n-d}\bigg\rfloor +e =e_0 +\bigg\lfloor\frac{2e_0-2}{d-2}\bigg\rfloor\eqno{(A.1)}$$
is the intersection of $[d-1, \infty)$ with a union of
%$N(d,n)$
at least $$\min\bigg\{\big(n-d-\frac 12\big)\big(d-3\big)- \frac{15}2, \;(n-d)(d-5)-2,\; (n-d+1)(d-5)\bigg\}$$ congruence classes$\mod d(n-2)$.
% where $(n-d)(d-2)\leq N(d,n)\leq (n-d)(d-1)+2.$
%the number of congruence classes  is $(n-d)(d-2)+X$.
%$$\begin{aligned} &(n-d)(d-2)+X \text{  for n-d even, or n-d odd and g odd,}\qquad\qquad\qquad\qquad\qquad\qquad\qquad\\
%& (n-d)(d-2)+2+X\text{  for n-d odd, g even, and d even,}\\
%& (n-d)(d-2)+3+X\text{  for n-d odd, g even, and d odd}.\end{aligned}$
}

\bigskip

\noindent{\it Remark 1.}  For $n, d >> 0$ the density of the accessible $e$ is about $(n-d)/n$.

\bigskip

\noindent{\it Remark 2.} The formula for the lower bound on $N(d,n)$, the number of congruence classes of $e\mod d(n-2)$ will be obtained by applying Facts 1-3 to count the
number of permissible $c\in [0,d-3]$  in display (A.3) below.
%Counting $c$'s is done through counting the number of congruence
%classes of $c \mod n-d+1$ assuming all classes have the same number of elements, and assuming $c>0$. This leads
%to the (nonnegative) `error' term $N(d,n)-(n-d)(d-\frac 52)$ and its upper bound in Remark 2 at the end of the proof of Theorem A.1. For the case $c=0$, see the discussion in Remark 4.
%The error terms $X_1, X_2, X_3$ occur. Because to count the number of permissible $c\in [0, d-2]$ we count the number of congruence classes of $c\mod n-d+1$ by assuming all classes have the same number of element. So there may be up to $n-d$ classes do not have the average number of elements.

\bigskip

\noindent{\bf Proof of Theorem A.1.}

\medskip

\noindent Fix $n, d$ as in Theorem A.1. Let
$$b=n-d+1.
\eqno{(A.2)}$$
Dividing $e_0$ by $d-2$, then dividing the quotient obtained by $b$, we have
$$e_0=k(d-2)b+r(d-2)+c, \eqno{(A.3)}$$
where
$$0\leq c\leq d-3,\;\text{ and } \; 0\leq r\leq b-1. \eqno{(A.4)}$$
Hence
$$\frac{2e_0-2}{d-2}=2kb+2r+\frac{2c-2}{d-2}. \eqno{(A.5)}$$

\bigskip

\noindent {\it Since we consider the lower bound on the number of $e$'s, we may assume that $c>0$.}

\bigskip

\noindent There are the following two cases of the integral part of $\frac{2e_0-2}{d-2}$.

\smallskip

\noindent
{\it Case} (a). $0<c<\frac d2$.

\noindent In this case we have $\frac{2c-2}{d-2}<1$. Therefore, from display (A.5)
$$\bigg\lfloor\frac{2e_0-2}{d-2}\bigg\rfloor=2kb+2r .\eqno{(A.6.a)}$$
{\it Case} (b). $c\geq\frac d2$.

\noindent In this case we have  $1\leq\frac{2c-2}{d-2}<2$. Therefore,
$$\bigg\lfloor\frac{2e_0-2}{d-2}\bigg\rfloor=2kb+2r+1,\eqno{(A.6.b)}$$
and the fractional part of $\frac{2e_0-2}{d-2}$ is
$$\bigg\{\frac{2e_0-2}{d-2}\bigg\}=\frac{2c-2-(d-2)}{d-2}.\eqno{(A.7)}$$

\bigskip

\noindent Coming back to display (A.1), we let $\varepsilon$ be the fractional part of $\frac{-de_0+e}{n-d}$,

\noindent i.e.,
$$\frac{-de_0+e}{n-d} = \bigg\lfloor \frac{-de_0+e}{n-d}\bigg\rfloor +\varepsilon.\eqno{(A.8)}$$
In particular, $\varepsilon<1$.

\bigskip

\noindent Putting displays $(A.1), (A.2), (A.3), (A.6), (A.8)$ together, we have
$$\text{Case (a).   }\quad e=d(n-2)k+rd+c+\frac {r(d^2-3d)+c(d-1)}{b} + \varepsilon\frac{b-1}{b} \qquad\qquad\qquad\eqno{(A.9)}$$
 $$
\text{Case (b).   }\quad e=d(n-2)k+rd+c+1+ \frac {r(d^2-3d)+c(d-1)-1}{b} + \varepsilon\frac{b-1}{b} \qquad\qquad\qquad$$
for Cases (a) and (b) respectively.

\bigskip

We want to count the values $e$ expressed in display (A.9) with all possible $(c, r) \in [1, d-3]\times[0, b-1]$ by counting congruence classes of $e$ mod $d(n-2)$.

\bigskip

We will give the argument for Case (a) only, since the argument for Case (b) is identical. Let $$E(c, r, \varepsilon)=rd+c+\frac {r(d^2-3d)+c(d-1)}{b} + \varepsilon\frac{b-1}{b}.$$

\noindent{\it Claim 1.} If $(c,r)\not= (c_1, r_1)$, then $E(c, r, \varepsilon)\not= E(c_1, r_1, \varepsilon_1)$ as real numbers.

\smallskip

\noindent{\it Proof of Claim 1.}

Assume $E(c, r, \varepsilon)= E(c_1, r_1, \varepsilon_1)$. Then
$$(r_1-r)\bigg(d+\frac{d^2-3d}{b}\bigg)=(c-c_1)\bigg(1+\frac{d-1}{b}\bigg) +(\varepsilon-\varepsilon_1)\frac{b-1}{b}.\eqno{(A.10)}$$We may assume $r_1-r\geq 1$. Hence the left-hand-side of display (A.10) gives
$$(r_1-r)\bigg(d+\frac{d^2-3d}{b}\bigg)\geq 1\cdot \bigg(d+\frac{d^2-3d}{b}\bigg),$$
while, by (4) and that $\varepsilon, \varepsilon_1 \in [0,1)$, the right-hand-side of display (A.10) gives
$$(c-c_1)\bigg(1+\frac{d-1}{b}\bigg) +(\varepsilon-\varepsilon_1)\frac{b-1}{b}\leq (d-3)\cdot \frac {b+d-1}{b} +\frac {b-1}{b}.$$ This is a contradiction.

\bigskip
\noindent{\it Claim 2.} $E(c, r, \varepsilon) < d(n-2)$.

\smallskip

\noindent This is clear, because again, by displays (A.4) and (A.8)
$$\begin{aligned}E(c, r, \varepsilon)&\leq(b-1)\bigg(d+\frac{d^2-3d}{b}\bigg)+(d-3)\bigg(1+\frac{d-1}{b}\bigg) +
\frac{b-1}{b}\\&< d(b+d-3)\\&=d(n-2).\end{aligned}$$
From Claim 1 and Claim 2, we conclude
$$\text{  if  } (c,r)\not= (c_1, r_1), \text{  then  } E(c, r, \varepsilon)\not\equiv E(c_1, r_1, \varepsilon_1)\mod d(n-2),\qquad\qquad\qquad$$i.e., for $e=e(c,r)$ in display (A.9)
$$\text{if  } (c,r)\not= (c_1, r_1), \text{  then  } e(c,r)\not\equiv e(c_1, r_1)\mod d(n-2).\qquad\eqno{(A.11)}$$

\bigskip

Next, we want to count the permissible $(c, r) \in [1, d-3]\times[0, b-1]$.
%\{0, 1, 2,\ldots, d-3\}\times \{0, 1, 2,\ldots, b-1\}$.
As before, we give the argument for Case (a).

\noindent Since $e\in \mathbb Z$, we may let
$$\frac {r(d^2-3d)+c(d-1)}{b} + \varepsilon\frac{b-1}{b}=m+1, \text{ for } m\in \mathbb Z.\eqno(A.12)$$
Hence
$$\begin{aligned} &m< \frac {r(d^2-3d)+c(d-1)}{b}\qquad\qquad\qquad\qquad\qquad\qquad\qquad\qquad\qquad\qquad\qquad\qquad\\ \text{and}\qquad\qquad&\qquad\qquad\qquad\qquad\qquad\\
&\varepsilon=\frac{(m+1)b-\big(r(d^2-3d)+c(d-1)\big)}{b-1}.\end{aligned}\eqno(A.13)$$
%Since $\varepsilon$ is the fractional part of $\frac{-de_0+e}{n-d}$,
By display (A.8), $\varepsilon<1$, which is equivalent to
$$m< \frac {r(d^2-3d)+c(d-1)-1}{b}.$$
So we want to rule out those $(c,r)$ such that
$$\frac {r(d^2-3d)+c(d-1)-1}{b}=m,$$i.e., we want to rule out $(c, r) \in [1, d-3]\times[0, b-1]$
%\{0, 1, 2,\ldots, d-3\}\times \{0, 1, 2,\ldots, b-1\}$
such that
$$r(d^2-3d)+c(d-1)\equiv 1\mod b.\eqno(A.14)$$
By display (A.2), solving the congruence equation (A.14) is the same as solving
$$cn\equiv -r(n+1)(n-2)+1 \mod b.\eqno(A.15.a)$$

\medskip

We will use the following facts about the
 congruence equation
$$ax\equiv d \mod b\eqno{+}$$
\noindent
{\it Fact 1.}  Equation ${+}$ is solvable if and only if g:=gcd$(a,b)$ divides $d$.

\smallskip

\noindent
{\it Fact 2.} Assume $g|d$, and let $b':=b/g$. If we consider the solution $x$ as an integer, then $x$ is unique in any interval of size $b'$.

\smallskip

\noindent
{\it Fact 3.} For $\mathcal{C}\geq b'$, the number of solutions of ${+}$ in $[1, \mathcal{C}]$ is $\big\lfloor\mathcal{C}/b'\big\rfloor$ or $\big\lfloor\mathcal{C}/b'\big\rfloor+1$.

\bigskip

Coming back to congruence equation (A.15.a), we let $g=$gcd$(n,b)$. Counting the numbers of $r\in[0,b-1]$ such that $$g|-r(n+1)(n-2)+1\eqno(A.16)$$  is the same as counting $r$ satisfying
$$ g| 2r+1\eqno(A.17)$$i.e., counting the number of $r$ such that
$$2r\in \{g-1, 2g-1,\ldots, 2b'g-1\},\text{ where }b'=\frac {b}{g}.\eqno(A.18)$$

\medskip

\noindent
{\bf I.a.} Assume $\mathcal{C}_a \geq b'$, where $\mathcal{C}_a= \big|[1, \frac d2)\bigcap\mathbb Z|$.

\medskip
\noindent{\it Case} (I.a.i.) $b$ is odd.

\noindent In this case, $g=$gcd$(n,b)$ is odd and $$\begin{aligned}&\{2r: r\in [0, b-1] \text{ and r satisfies (A.15.a) }\} \qquad\\=\;&\{g-1, 3g-1,\ldots, (2b'-1)g-1\}.\end{aligned}\eqno(A.19.a)$$

\noindent {\it There are at most $\big(\big\lfloor\mathcal{C}_a /b'\big\rfloor+1\big)\;b'\leq \mathcal{C}_a+b'$ pairs of $(c,r)$.}

\bigskip
\noindent{\it Case}(I.a.ii.1.) $b$ is even and $g$ is even. (Hence $n$ is even.)

\medskip
\noindent  {\it There is no $r$ satisfying display (A.17).}

\bigskip
\noindent{\it Case} (I.a.ii.2.) $b$ is even and $g$ is odd. (Hence $n$ is odd.)
$$\begin{aligned}&\{2r: r\in [0, b-1] \text{ and r satisfies (A.15.a) }\} \qquad\qquad\qquad\\=\;&\{g-1, 3g-1,\ldots,(2b'-1)g-1\}.\end{aligned}\eqno(A.20.a)$$
 {\it There are at most $\big(\big\lfloor\mathcal{C}_a /b'\big\rfloor+1\big)\;b'\leq \mathcal{C}_a+b'$ pairs of $(c,r)$.}

\bigskip

\noindent
{\bf I.b.} Assume $\mathcal{C}_b \geq b'$, where $\mathcal{C}_b= \big|[\frac d2, d-3]\bigcap\mathbb Z|$.
For Case (b), we have
$$cn\equiv -r(n+1)(n-2)+2 \mod b,\eqno(A.15.b)$$and hence the following cases.

\medskip

\noindent{\it Case} (I.b.i.) $b$ is odd. (Hence $g$ is odd.) $$\begin{aligned}&\{2r: r\in [0, b-1] \text{ and r satisfies (A.15.b) }\} \qquad\\=\;&\{2g-2, 4g-2,\ldots, 2b'g-2\}.\end{aligned}\eqno(A.19.b)$$
  \noindent{\it There are at most $\big(\big\lfloor\mathcal{C}_b /b'\big\rfloor+1\big)\;b'\leq \mathcal{C}_b+b'$ pairs of $(c,r)$.}

\bigskip
\noindent{\it Case} (I.b.ii.1.) $b$ is even and $g$ is even.
$$\begin{aligned}&\{2r: r\in [0, b-1] \text{ and r satisfies (A.15.b) }\} \qquad\\=\;&\{g-2, 2g-2,\ldots,2b'g-2\}.\end{aligned}\eqno(A.21.b)$$
\noindent {\it There are at most $\big(\big\lfloor\mathcal{C}_b /b'\big\rfloor+1\big)\;2b'\leq 2\mathcal{C}_b+2b'$ pairs of $(c,r)$.}

\bigskip
\noindent{\it Case} (I.b.ii.2.) $b$ is even and is $g$ odd.
$$\begin{aligned}&\{2r: r\in [0, b-1] \text{ and r satisfies (A.15.b) }\} \qquad\\=\;&\{2g-2, 4g-2,\ldots,2b'g-2\}.\end{aligned}\eqno(A.20.b)$$
\noindent {\it There are at most $\big(\big\lfloor\mathcal{C}_b /b'\big\rfloor+1\big)\;b'\leq \mathcal{C}_b+b'$ pairs of $(c,r)$.}

\bigskip

\noindent
{\bf II.a.} {\it Assume $\mathcal{C}_a < b'.$}

\medskip

\noindent For each $r$, there is at most one solution $c$ in $[1, \mathcal{C}_a].$

\noindent Hence

\noindent
{\it Case} (II.a.i.) $b$ is odd. {\it There are at most $b'$ pairs of $(c,r)$.}

\smallskip
\noindent{\it Case} (II.a.ii.1.) $b$ is even and $g$ is even. {\it There is no $r$ satisfying display (A.17).}

\smallskip
\noindent{\it Case} (II.a.ii.2.) $b$ is even and $g$ is odd.  {\it There are at most $b'$ pairs of $(c,r)$.}

\bigskip

\noindent
{\bf II.b.} {\it Assume $\mathcal{C}_b < b'.$ }

\smallskip
\noindent{\it Case} (II.b.i.) $b$ is odd. {\it There are at most $b'$ pairs of $(c,r)$.}

\smallskip
\noindent{\it Case} (II.b.ii.1.) $b$ is even and $g$ is even. {\it There are at most $2b'$ pairs of $(c,r)$.}

\smallskip
\noindent{\it Case} (II.b.ii.2.) $b$ is even and $g$ is odd. {\it There are at most $b'$ pairs of $(c,r)$.}

 \bigskip

\bigskip
Summing up Cases (a) and (b), and using the facts that
$$\begin{aligned}&(1). \;\mathcal{C}_a+\mathcal{C}_b= d-3,\qquad\qquad\qquad\qquad\qquad\qquad\qquad\qquad\qquad\qquad\qquad\qquad\\
&(2).\; \mathcal{C}_a=\frac d2-1, \text{  if  d is even},\\
&\qquad\mathcal{C}_a=\frac {d-1}2, \text{  if d is odd},\\&\qquad\mathcal{C}_b=\frac d2-2, \text{  if  d is even},\\
&\qquad\mathcal{C}_b=\frac {d+1}2-3, \text{  if d is odd}.\end{aligned}$$
 Taking off the bad pair $(c,r)$ from Cases (a) and (b), we have that the number of the permissible $(c,r) \in [1, d-3]\times [0, b-1]$ is at least
 $$\begin{aligned}&(I_a\,\&\,I_b).\text{When $\mathcal{C}_a>\mathcal{C}_b\geq b'$},\\
 &\quad\quad\quad\;\;\,(1). (b-1)(d-5)-2, \text{ if b is odd, or b is even and g is odd,}\,\qquad\qquad\qquad\qquad\qquad\qquad\qquad\qquad\qquad\qquad\qquad\qquad\\
  &\quad\quad\quad\;\;\,(2). (b-1)(d-5)-1, \text{ if b is even, g is even, and d is even, }\,\qquad\qquad\qquad\qquad\qquad\qquad\qquad\qquad\qquad\\
  &\quad\quad\quad\;\;\,(3). (b-1)(d-5), \;\;\;\;\;\;\text{ if b is even, g is even, and d is odd. }\end{aligned}$$
  $$\begin{aligned}&({I}_a\,\&\, {II}_b).\text{When $\mathcal{C}_a\geq b'>\mathcal{C}_b$},\\
 &\quad\quad\quad\;\;\,(1). (b-\frac 32)(d-3)-\frac{15}2, \text{ if b is odd, or b is even and g is odd,}\,\qquad\qquad\qquad\qquad\qquad\qquad\qquad\qquad\qquad\qquad\qquad\qquad\\
  &\quad\quad\quad\;\;\,(2). (b-1)(d-3)-1, \;\;\;\text{ if b is even, g is even, and d is even, }\,\qquad\qquad\qquad\qquad\qquad\qquad\qquad\qquad\qquad\\
  &\quad\quad\quad\;\;\,(3). (b-1)(d--3)-2, \;\text{ if b is even, g is even, and d is odd. }\end{aligned}$$
  $$\begin{aligned}(II_a\,\&\,II_b).&\text{When $b'>\mathcal{C}_a>\mathcal{C}_b,$ the number of the permissible pairs is at least}\quad\\&\;
  b(d-5).\qquad\qquad\qquad\qquad\qquad\qquad\qquad\qquad\qquad\qquad\qquad\qquad \end{aligned}$$

 \noindent Combining the above, we conclude
  the proof of Theorem A.1.$\qed$

\bigskip

\noindent {\it Remark 3.} The estimates can be improved by $b=n-d+1$, if $\mathcal{C}_a$ or $\mathcal{C}_b$ is a multiple of $b'=\frac b{\text{gcd}(b,n)}$. For example, in $I_a\,\&\,I_b$ (when $\mathcal{C}_a>\mathcal{C}_b\geq b'$), the number of permissible pairs $(c,r)\in [1, d-3]\times [0, b-1]$ is at least
$$\begin{aligned} &(b-1)(d-4) -1\text{  for b odd, or b even and g odd,}\qquad\qquad\qquad\qquad\qquad\qquad\qquad\\
& (b-1)(d-3)+1\text{  for b even, g even, and d even,}\\
& (b-1)(d-3)+2\text{  for b even, g even, and d odd}.\end{aligned}$$
\bigskip

\noindent {\it Remark 4.} Suppose $b=n+1-d\not= 2$. If $c=1$, then at least half of $r\in[0,b-1]$ are permissible. This can be seen in equation (A.15.a) and Facts 1 and 3. Fact 1 implies that $n-2$ and $b$ are relatively prime. If all $r\in[0,b-1]$ satisfy (A.15.a), then let $r=0$ and $r=1$, we have $b=2$. Together with Fact 3, we see that the number of solutions $r$ is a proper factor of $b$.

\bigskip
\noindent {\it Remark 5.} In $({I}_a\& {II}_b)$ the condition $\mathcal{C}_a\geq b'>\mathcal{C}_b$ implies that $$\begin{aligned}&\frac {d-3}2 =b'-\frac 12, \; \qquad\;\;\,\text{ if } d \text{ even},\\& \frac {d-3}2 =b'-1, \text{ or } b', \; \;\;\text{ if }d \text{  odd}.
\end{aligned}$$
\bigskip

\noindent{\bf Theorem A.2.} (i). {\it For all integers $n, d,\;$ $3<d<n-1$, set $$q(e)=\bigg\lfloor \frac{e(n+1-d)-2}{n-2}\bigg\rfloor+1.$$There are at least $\frac {n+1-d}2$ ray congruence 
	classes (arithmetic progressions) $\mod d(n-2)$ of $e$ which are both accessible and point-minimal, i.e. beside condition $(A.1)$, $e$ also satisfies
$$  q(e-1) < q(e).\eqno{(A.22)} $$}
\noindent (ii).
 {\it The set of values $q(e)$ with $e$ accessible and point-minimal contains
$(n+1-d)/2$ many ray congruence classes$\mod (n+1-d)d$.}

\medskip

\noindent
(iii). {\it If integers $n, d,\;$ $ 3<d<n-1$ such that
 $$\rm{gcd}\big((n+1)(n-2),\; n+1-d\big)=1, \eqno{(A.23)}$$then there are at least $f(n,d)$ ray congruence classes$\mod d(n-2)$  (respectively,$\mod d(n+1-d)$) of accessible and point-minimal $e$ (resp. of values q(e) over those),
where
$$\begin{aligned}&f(n,d)\\=&\min\bigg\{\frac{(n-d)(n-d-1)}2 -2,  \frac{(d+1)(d-6)}2, \frac{(n-d+1)(n-d-2)}2 -3\bigg\}.   \end{aligned}\eqno{(A.24)}$$}
%is the intersection of $[d-1, \infty)$ with a union of
%$N(d,n)$
%at least $$\min\bigg\{\big(n-d-\frac 12\big)\big(d-3\big)- \frac{15}2, \;(n-d)(d-5)-2\bigg\}$$ congruence classes$\mod d(n-2)$.}

\bigskip
\noindent {\it Remark 6.} By definition, the point-minimal condition does not affect the set of $q$-values. Thus,
$$\{q(e): e \text{ is accessible and point-minimal} \}=\{q(e): e \text{ is point-minimal}\}.$$

\bigskip
\noindent {\it Remark 7.} For $n, d \gg 0$ the density of the accessible and point-minimal $e$ is about $\frac{n-d}{2nd}$ in Theorem A.2(i) or $\min \big(\frac{(n-d)^2}{2nd}, \frac{d}{2n}\big)$ in Theorem A2(iii), and the density of the set of $q$ values is $\frac 1{2d}$ for (ii) and $\min\big(\frac{n-d}{2d}, \frac d{2(n-d)}\big)$ for (iii).

\bigskip

\noindent{\bf Proof of Theorem A.2.}

\medskip

\noindent
We will use the set up in the proof of Theorem A.1. By the expression of $e$ in display (A.9), inequality (A.22) is equivalent to the following inequality

\bigskip

\noindent {\it Case} (a). $0<c<\frac d2$.
$$ dbk+rd+\bigg\lfloor\frac{cn+\varepsilon(b-1)-2-b}{n-2}\bigg\rfloor < dbk+rd+\bigg\lfloor\frac{cn+\varepsilon(b-1)-2}{n-2}\bigg\rfloor  \qquad\qquad\qquad$$
\noindent{\it Case} (b). $c\geq\frac d2$.
$$ dbk+rd+\bigg\lfloor\frac{cn+\varepsilon(b-1)-3}{n-2}\bigg\rfloor < dbk+rd+\bigg\lfloor\frac{cn+b+\varepsilon(b-1)-3}{n-2}\bigg\rfloor  \qquad\qquad\qquad $$
Since there is no restriction on $k$, as long as one pair of $(c,r)$ such that inequality holds, there are infinitely many pair. Hence we want to find $(c,r)$ such that
$$0<c<\frac d2, \;\; \text{  and  } \;\;\bigg\lfloor\frac{cn+\varepsilon(b-1)-2-b}{n-2}\bigg\rfloor < \bigg\lfloor\frac{cn+\varepsilon(b-1)-2}{n-2}\bigg\rfloor,  \qquad\qquad\qquad\eqno{(A.25.a)}$$
or
$$c\geq\frac d2, \;\;\text{  and  }\;\;\bigg\lfloor\frac{cn+\varepsilon(b-1)-3}{n-2}\bigg\rfloor < \bigg\lfloor\frac{cn+b+\varepsilon(b-1)-3}{n-2}\bigg\rfloor.  \qquad\qquad\qquad \eqno{(A.25.b)}$$
We will only prove the theorem for Case (a), since Case (b) is similar.

\medskip

In display (A.25.a), it is straightforward that
$$L:=\frac{cn+\varepsilon(b-1)-2-b}{n-2}>c-1, \;\text{ and } \;R:=\frac{cn+\varepsilon(b-1)-2}{n-2}>c.$$
The condition $0<c<\frac d2$ implies
$$R=\frac{cn+\varepsilon(b-1)-2}{n-2}<c+1.$$
Hence
$$\lfloor R \rfloor=\bigg\lfloor\frac{cn+\varepsilon(b-1)-2}{n-2}\bigg\rfloor =c, \qquad\qquad\qquad\eqno{(A.26.a)}$$
and the permissible pairs $(c,r)$ for $0<c<\frac d2$ are exactly those satisfy
$$L=\frac{cn+\varepsilon(b-1)-2-b}{n-2}<c,$$
i.e.
$$\varepsilon < \frac{b+2-2c}{b-1}.\qquad\qquad\qquad\eqno{(A.27.a)}$$It is clear that (A.27.a) holds for $c=1$, and for any $r\in [0,b-1]$. So does (A.22). On the other hand, by Remark 4, (A.1) holds for at least half of $r\in [0,b-1]$. Hence,  parts (i) and (ii) are  proved, and we may assume $c\in [2, \frac d2)$ for the proof of part (iii).

\smallskip

\noindent By the expression of $\varepsilon$ in (A.13), display (A.27.a) is equivalent to
$$(m+1)b-r(d^2-3d)-c(d-1) < b+2-2c,\qquad\qquad\qquad\eqno{(A.28.a)}$$
where $m$ is as in (A.12) and (A.13). So we want to rule out $(c,r)$ such that
$$ \frac {r(d^2-3d)+c(d-3)+2}{b}\leq m< \frac {r(d^2-3d)+c(d-1)}{b},\eqno(A.29.a)$$
i.e. we want to rule out $(c,r)$ such that
$$r(d^2-3d)+c(d-1)\equiv K\mod b, \;\text{ for } K=2,\ldots, 2c-2,$$
(The case $K=1$ was done in Theorem A.1.)

\noindent or equivalently
$$r(n+1)(n-2)\equiv K-cn\mod b, \;\text{ for } K=2,\ldots, 2c-2.\eqno(A.30.a)$$
For each $c$, there are $2c-3$ choices of $K$, and for each $K$, the solution $r$ is unique under assumption (A.23). So the total number of $(c,r)\in [2, \mathcal C_a]\times [0,b-1]$ to rule out is
$$\begin{aligned}\quad&1+3+\cdots+(b-2)+(b-1)\big(\mathcal{ C}_a -\frac{b+1}2\big)\\= &(b-1)\mathcal{ C}_a-\frac{(b-1)(b+3)}{4},\;\;\text{  if $b$ is odd, and }\mathcal{C}_a\geq \frac{b+1}2, \qquad\qquad\qquad\qquad\qquad\qquad\qquad\qquad\qquad\qquad\qquad\qquad\\\text{and}\quad&\qquad\qquad\qquad
\\
\quad\qquad&1+3+\cdots+(2\mathcal{C}_a-3)\\= &(\mathcal{C}_a-1)^2, \qquad\qquad\qquad\qquad\text{  if }\mathcal{C}_a< \frac{b+1}2. \qquad\qquad\qquad\qquad\qquad\qquad\qquad\qquad\qquad\qquad\qquad\qquad
\end{aligned}$$

\noindent
For Case (b), $c\geq\frac d2$, we have the corresponding displays
$$\lfloor R_b \rfloor=\bigg\lfloor\frac{cn+b+\varepsilon(b-1)-3}{n-2}\bigg\rfloor =c, \qquad\qquad\qquad\eqno{(A.26.b)}$$
$$\varepsilon < \frac{b+d-2c}{b-1}.\qquad\qquad\qquad\eqno{(A.27.b)}$$
$$(m+1)b-r(d^2-3d)-c(d-1) < b+d-2c,\qquad\qquad\qquad\eqno{(A.28.b)}$$
$$ \frac {r(d^2-3d)+c(d-3)-2c+d-1}{b}\leq m< \frac {r(d^2-3d)+c(d-1)-1}{b},\eqno(A.29.b)$$
$$r(n+1)(n-2)\equiv K-cn\mod b, \;\text{ for } K=3,\ldots, 2c-d+1.\eqno(A.30.b)$$

 \noindent Summing up Case (a) and Case (b), the number of pairs $(c,r)\in [2,d-3]\times [0, b-1]$ to rule out is at most
 $$\begin{aligned}&(1_a\,\&\,1_b).\text{When $\mathcal{C}_a>\mathcal{C}_b\geq \frac {b+1}2$},\\
 &\quad\quad\quad\;\;\,(1). (b-1)(d-3)-\frac {(b-1)(b+3)}2, \text{ for b  odd and d  even,}\,\qquad\qquad\qquad\qquad\qquad\qquad\qquad\qquad\qquad\qquad\qquad\qquad\\
  &\quad\quad\quad\;\;\,(2). (b-1)(d-3)-\frac {(b-1)(b+2)}2, \text{ for  d odd, }\,\qquad\qquad\qquad\qquad\qquad\qquad\qquad\qquad\qquad\\
  &\quad\quad\quad\;\;\,(3). (b-1)(d-3)-\frac {b^2+2b-4}2, \;\;\;\text{ for b and d  even . }\end{aligned}$$
  $$\begin{aligned}&({1}_a\,\&\, {2}_b).\text{When $\mathcal{C}_a\geq \frac{b+1}2>\mathcal{C}_b$},\\
 &\quad\quad\quad\;\;\,(1). \frac {b^2-4b+5}2, \text{ when  } b=d-3,\,\qquad\qquad\qquad\qquad\qquad\qquad\qquad\qquad\qquad\qquad\qquad\qquad\\
  &\quad\quad\quad\;\;\,(2). \frac {b^2-5b+8}2, \text{ when  } b=d-2.\,\qquad\qquad\qquad\qquad\qquad\qquad\qquad\qquad\qquad\end{aligned}$$
  $$\begin{aligned}&({2}_a\,\&\, {2}_b).\text{When $\frac{b}2>\mathcal{C}_a>\mathcal{C}_b$},\\
 &\quad\quad\quad\;\;\,(1).  \frac {d^2-10d+26}2, \text{ if d is even, }\,\qquad\qquad\qquad\qquad\qquad\qquad\qquad\qquad\qquad\\
  &\quad\quad\quad\;\;\,(2). \frac {d^2-9d+22}2,\text{ if d is odd. } \end{aligned}$$
   
  \noindent Note that from displays (A.25.a), (A.25.b), (A.26.a), and (A.26.b), the number of values $q(e)\mod db$ is the same as the number of permissible pair $(c,r)$.
  
  \noindent In each case of $(1_a\,\&\,1_b), (1_a\,\&\,2_b)$, and $(2_a\,\&\,2_b)$ above, we take the maximum and subtract it from the minimum in Theorem A.1. Applying Remark 5, we prove Part (iii). $\qed$

\bigskip

\end{document}